\documentclass[11pt]{article}
\usepackage[T1]{fontenc}
\usepackage[utf8]{inputenc}
\usepackage[english]{babel}

\usepackage{amsmath, amsthm, amssymb, amsfonts, enumerate,mathrsfs,cases,bbm}
\usepackage[colorlinks=true,linkcolor=blue,urlcolor=blue]{hyperref}
\usepackage{dsfont}
\usepackage{color}
\usepackage{geometry}
\usepackage{todonotes}
\usepackage{listings}
\usepackage[normalem]{ulem}
\usepackage{hyperref}
\hypersetup{
    colorlinks,
    citecolor=blue,
    filecolor=blue,
    linkcolor=blue,
    urlcolor=blue
}
\usepackage[round]{natbib}

\newtheorem{theorem}{Theorem}[section]
\newtheorem{remark}[theorem]{Remark}

\newtheorem{proposition}[theorem]{Proposition}

\newtheorem{definition}[theorem]{Definition}
\newtheorem*{notation*}{Notation}  

\def \E{\mathbb{E}}
\def \V{\mathbb{V}}

\DefineNamedColor{named}{cgreen}{rgb}{0.0, 0.42, 0.24}

\DeclareMathOperator*{\argmax}{arg\,max}

\numberwithin{equation}{section}

\begin{document}
\author{René Aïd\thanks{Universit\'e Paris-Dauphine, PSL Research University, rene.aid@dauphine.psl.eu}\quad
Ofelia Bonesini\thanks{Universit\`a  degli Studi di Padova, Dipartimento di Matematica ``Tullio Levi-Civita'', bonesini@math.unipd.it, ORCID 0000-0001-9294-6079.}\quad
Giorgia Callegaro\thanks{Universit\`a degli Studi di  Padova, Dipartimento di Matematica ``Tullio Levi-Civita'', gcallega@math.unipd.it, ORCID 0000-0001-9026-5261.}\quad
Luciano Campi\thanks{Universit\`a degli Studi di Milano, Dipartimento di Matematica ``Federigo Enriques'', luciano.campi@unimi.it, ORCID 0000-0002-3956-6795.}
}
\title{A McKean-Vlasov game  \\ of commodity production, consumption and trading}

\date{\today}
\maketitle

\begin{abstract} 
We propose a model where a producer and a consumer can affect the price dynamics of some commodity controlling drift and volatility of, respectively, the production rate and the consumption rate. 
We assume that the producer has a short position in a forward contract on $\lambda$ units of the underlying  at a fixed price  $F$, while the consumer has the corresponding long position.  Moreover, both players are risk-averse with respect to their financial position and their risk aversions are modelled through an integrated-variance penalization. 
We study the impact of risk aversion on the interaction between the producer and the consumer as well as on the derivative price.
In mathematical terms, we are dealing with a two-player linear-quadratic McKean-Vlasov stochastic differential game. Using methods based on the martingale optimality principle and BSDEs, we find a Nash equilibrium and characterize the corresponding strategies and payoffs in semi-explicit form.
Furthermore, we compute the two  indifference prices (one for the producer and one for the consumer) induced by that equilibrium and we determine the quantity $\lambda$ such that the players agree on the price.
Finally, we illustrate our results with some numerics. In particular, we focus on how the risk aversions and the volatility control costs of the players affect the derivative price.
\end{abstract}


\smallskip

{\textbf{MSC2020 subject classification}}: 
 49N10,
 91A15, 
 91G30.

\section{Introduction}\label{intro}

In this paper, we develop an economic model of a commodity market where a representative producer interacts with a representative processor who buys the commodity and transforms it into a final product sold to the retail market (think of crude oil production transformed into gasoline or wheat transformed into bread). For the sake of simplicity, the processor will be referred to as consumer from now on. In our model, the production and the consumption rates are described as It\^o processes driven each by an independent Brownian motion and whose coefficients are controlled by, respectively, the producer and the consumer. We stress that in our model the producer can control, in particular, the volatility of the production rate (by investing in devices making the production more reliable), and similarly the consumer can control the one of the consumption rate (by investing in storage devices, for instance). Further, the players are risk-averse (see below for details) and they are linked by a financial derivative in the commodity, a plain forward agreement on price and volume exchanged. For some motivations on the control of volatility, we refer the reader to the paper by \cite{ACC20}, which focuses on the interaction between a producer controlling the drift of the spot price and a trader controlling the volatility, and exchanging a quadratic derivative. In that paper, it was shown that when the trader is short in the derivative, he would increase the volatility of the spot price in order to get a higher price of the derivative sold to the producer. In the present setting, we are interested in the joint effect of the costs of controlling the volatility of production or consumption rates and the players' risk aversion parameters on the ``agreement indifference price''.
Indeed, when only one player has market power, the effect of the parameters on the forward price is clear. On the other hand, when the two players interact, the joint effect is not obvious. In this paper, we are interested in the outcome of the combined effect on the forward price of the relative risk aversions and the volatility control costs of the producer and the consumer.

Both players have market power on the spot price of the commodity: the spot price depends linearly on production and consumption rates so that the higher the rate of production, the lower the spot price and the higher the rate of consumption, the higher the spot price. Furthermore, they agree to exchange a forward contract with finite maturity $T$ over a certain quantity $\lambda$ of the commodity that will be determined at equilibrium together with its price $F$. This setting is inspired from the seminal papers of  \cite{A92} and \cite{AV93}, where the authors establish the mitigating effect of forward agreement on the exercise of producers market power.

In our framework, since production and consumption rates are driven by two independent Brownian motions and there is only one tradable risky asset, i.e. the commodity spot price, the market is incomplete. Therefore, we define the forward price in the spirit of the indifference pricing approach (see the paper \cite{HH09} for an overview and \cite{BCK08} for an application to power markets). The players' goal is to maximize their respective objective functionals, which are expectations of the following main components: the profit from selling,  the sourcing costs (only for the consumer), the costs from exerting the controls, the forward contract payoff and, finally, the integrated variance of the market price of the derivative. 

The latter component describes the risk aversion both players have towards their financial position. More precisely, in this context where the agents can control the volatility of their state variable, the modelling of their risk aversion using utility functions (e.g. exponential utility) would lead to nonlinear PDEs which are difficult to handle. Hence, for technical convenience we turn to a sort of dynamic mean-variance criterion leading to the objective functionals described above.  
Mathematically speaking we are dealing with a two-player stochastic differential game with objective functionals of McKean-Vlasov type, i.e. depending on the laws of the state variables. 
Economically speaking, it means that both players act as speculators on the forward market, as they disconnect their forward position from their production or transformation profit. Although this feature of our model originates from a computational limitation induced by the linear-quadratic McKean-Vlasov game setting, there exists some evidence, documented by a stream of the economic literature, that large commodity players can act as speculators on their markets (see \cite{CX14} for such evidence and references on the subject of financiarisation of commodity markets).

This modeling approach for the risk aversion has been already investigated and used for portfolio selection by \cite{ZL00} and more recently by \cite{IP19} and \cite{LLP20}. Moreover, due to the fast development of mean-field games as a new framework to study stochastic differential games for a large number of players since the seminal papers by \cite{LL06,LL06a,LL07} and \cite{HMC06} (see also \cite{CHM17} for a survey), there has been a regain of interest for control problems of McKean-Vlasov dynamics. The latter, also known as mean field control, corresponds in some way to the limit of a sequence of stochastic control problems for a regulator willing to optimize the average expected payoff of a group of agents interacting through the empirical distribution of their states (see \cite{L17} and the two-volume book \cite{CD18}). In particular, the linear-quadratic case has been treated in \cite{G16}, \cite{BSYY16} and \cite{BP19}. 
Recently, stochastic differential games with both state dynamics and objective functionals of McKean-Vlasov type has been addressed in, e.g., \cite{MP18}, \cite{Cosso} and also \cite{FH20} for a Stackelberg game arising from an optimal portfolio liquidation problem. Although a large number of applications in economics and finance have been developed with mean field games and mean field control, the applications of games with finitely many players and McKean-Vlasov dynamics and objective functionals in economics is much more recent, hence less developed (see, e.g., \cite{ABP20}).\smallskip

We will analyse the model along the following program: first we will find a Nash equilibria for a fixed quantity $\lambda$ of the commodity exchanged through the forward contract with fixed price $F$; second, we will compute the indifference prices of the forward contract for the two players separately (they are going to depend on $\lambda$); third, we will compute the quantity $\lambda$ such that the two prices are equal, hence making the exchange compatible with the equilibrium found in the first step. This price will be called {\em agreement indifference price}.

This framework makes it possible to analyse the formation of the risk premium defined as the difference between the (unitary) agreement indifference price and the expected spot price of the commodity. The question of the determinants of the risk premium on commodity markets goes back (at least) to  Keynes's {\em Treatise on Money}, (1930). Keynes formulated the {\em normal backwardation theory}, i.e. the claim that forward prices should be lower than expected spot prices because risk-averse producers are willing to sell forward at a premium to avoid price risk. Presently, the {\em hedging pressure theory}  (see \cite{DNV00}, \cite{H88}, \cite{H88a}, \cite{H90}) provides explanation of the sign of the risk premium depending on the relative size of population types in the market (producers, storers, speculators) and their risk-aversion (see  \cite{ELV19} for a complete equilibrium model with mean-variance utility players explaining the different possible sign of the premia).

\paragraph{Mathematical results.} The main mathematical contribution of the paper (ref. Theorem \ref{thm_Nash}) consists in a complete description of a Nash equilibrium in open loop strategies of the two-player stochastic differential game arising from the interaction model described above. More in detail, we adopt the following resolution approach: first, we prove a suitable version of a verification theorem exploiting the weak martingale optimality principle; second, the verification theorem and the linear-quadratic structure of the game allows to provide a semi-explicit form for the best response map; third, a Nash equilibrium is found as a fixed point of the best response map with closed-form expressions for the equilibrium strategies and payoffs of both players up to solving numerically a Riccati system of ODEs. Once we have a Nash equilibrium at our disposal, computing the corresponding agreement indifference price together with the exchanged quantity at equilibrium is a pretty straightforward task. 

\paragraph{Economic insights.} First, we find that the forward agreement indifference price is higher (resp. lower) than the expected spot price when the producer is more (resp. less) risk-averse than the consumer. Because in our model, the players act as speculators on the forward market, a seller requires a higher forward price to enter in the agreement and a buyer asks for a lower price. The presence of market power of both players allows for the formation of an equilibrium. In that sense, our model is consistent with the economic intuition of the hedging pressure theory applied to a market populated with producers and consumers acting as speculators. Second,  we observe that producers can achieve the same agreement indifference price and the same trading volume either by having  high risk aversion and a low volatility  control  cost, or a low risk aversion and a high volatility  control cost. This effect manifests itself whatever the relative risk aversion of the producer and the consumer or the relative costs of volatility control. Nevertheless, it is more apparent when the volatility control costs are low. Thus, to the list of determinants of the sign of the risk premium of forward commodity price, one could add the costs of reducing the production uncertainty. For commodity where storage is utmost costly like electricity, reducing production uncertainty is highly costly and thus, leads to higher risk premium.

\paragraph{Organization of the paper.}  The paper is organised in the following way. The model is described in Section~\ref{sec:model} together with the definition of a forward agreement indifference price and quantity in Section \ref{sect_indiff_price}. The main result on the existence of a Nash equilibrium is given in Section~\ref{sec:nash}. The proof of the main result is given is Section~\ref{sec:proof}. Numerical results on the comparative static of the risk premium and the joint effect of risk aversion and volatility control costs are given in Section~\ref{sec:numeric}. 
 
\paragraph{Notations.}
	We denote by $\mathbb{R}_+$ (respectively $\mathbb{R}_-$) the closed semi-interval $[0, + \infty)$ (respectively $(-\infty, 0]$). 
	Given a function $f:\mathbb{R} \to S$, with $S$ a regular space, we denote its first derivative by $f'.$
	The expected value of a random variable $X$ will be equivalently denoted by $\mathbb{E}[X]$, as usual, or by $\bar{X}$, for brevity. 
	Let $(\Omega,\mathcal{F},\mathbb{P})$ be a probability space. Given a positive integer $d$, a strictly positive time horizon $T$  and a filtration $\mathbb{F}:=(\mathcal{F}_t)_{t \in [0,T]}$, we set
\begin{displaymath}
	\begin{split}
	&L^2([0,T], \mathbb{R}^d)
	:=\left\{ \varphi: [0,T] \to \mathbb{R}^d,\text{ s.t. } \varphi \text{ is measurable and }\int_0^T|\varphi_t|^2 dt < \infty   \right\}, \\&
	L^{\infty}([0,T], \mathbb{R}^d)
	:=\left\{ \varphi: [0,T] \to \mathbb{R}^d,\text{ s.t. } \varphi \text{ is measurable and }\sup_{t \in [0,T]}|\varphi_t| < \infty  \right\}, \\&
	L^2_{\mathcal{F}_T}(\Omega, \mathbb{R}^d)
	:=\left\{ \psi: \Omega \to \mathbb{R}^d,\text{ s.t. } \psi \text{ is } \mathcal{F}_T\text{-measurable and }\mathbb{E}\left[|\psi|^2 \right] < \infty   \right\},\\&
	L^2_{\mathbb{F}}(\Omega \times [0,T], \mathbb{R}^d)
	:=\left\{ \eta: \Omega \times [0,T] \to \mathbb{R}^d,\text{ s.t. } \eta \text{ is } \mathbb{F}\text{-adapted and }\mathbb{E}\left[\int_0^T|\eta_t|^2 dt \right] < \infty   \right\},\\&
	S^2_{\mathbb{F}}(\Omega \times [0,T], \mathbb{R}^d)
	:=\left\{ \eta: \Omega \times [0,T] \to \mathbb{R}^d,\text{ s.t. } \eta \text{ is } \mathbb{F}\text{-adapted and }\mathbb{E}\left[\sup_{t \in [0,T]}|\eta_t|^2 \right] < \infty   \right\}.
	\end{split}
\end{displaymath}

\section{The model}
We consider a stochastic game between a representative producer and a representative consumer. While the producer produces a good at a certain rate, the consumer buys the commodity and transforms it into a final good sold in the retail market.

\subsection{Market model} \label{sec:model}
We consider a finite time window $[0,T] $ and a probability space $(\Omega, \mathcal{F}, \mathbb{P}$) endowed with a two-dimensional  Brownian motion $(W,B)=\{(W_t,B_t)\}_{t \in [0,T]}$ and its natural filtration $\mathbb{F}=(\mathcal{F}_t)_{t \in [0,T]}$ augmented with the $\mathbb P$-null sets in $\mathcal F$.  The \emph{production rate} of the producer $\{q_t\}_{t \in [0,T]}$ evolves according to a dynamics given by
\[	
	dq_t  = u_t dt+   z_t  dW_t, \quad q_0>0,
\]
where $\{u_t\}_{t \in [0,T]}$ and  $\{z_t\}_{t \in [0,T]}$ are the  producer's strategies. The associated instantaneous costs are $\frac{k_p}{2}u_t^2$ and $\frac{\ell_p}{2}(z_t-\sigma_p)^2$, respectively,  with $k_p, \ell_p \geq 0$ and where $\sigma_p >0$ represents the nominal uncertainty in production without dedicated effort of the producer to reduce it. 
In a similar way, the \emph{consumption rate} (or selling rate to the retail market) of the consumer, $\{c_t\}_{t \in [0,T]}$, has dynamics given by
\[
	dc_t = v_t dt + y_t dB_t, \quad c_0>0.
\]
Here, $\{v_t\}_{t \in [0,T]}$ and  $\{y_t\}_{t \in [0,T]}$ are the \emph{consumer's strategies}, and the  associated instantaneous costs are, respectively, $\frac{k_c}{2}v_t^2$ and $\frac{\ell_c}{2}(y_t-\sigma_c)^2$, with $k_c, \ell_c \geq 0$ and $\sigma_c >0$. We assume a linear impact on the \emph{observed market price}, $\{S_t\}_{t \in [0,T]}$, namely $\{S_t\}_{t \in [0,T]}$  evolves according to
\[
	S_t := s_0 - \rho_p q_t +\gamma \rho_c c_t, \qquad s_0>0 
\]
with $\rho_p, \rho_c >0$ and $\gamma >0$ (the role of $\gamma$ will be clear in a few lines).
The instantaneous profits at time $t$ of the producer $\pi^p_t$ and of the consumer $\pi^c_t$  are  given by:	
		\begin{align*}
			\pi^p_t 
			:&= q_t S_t - \frac{k_p}{2}u_t^2 - \frac{\ell_p}{2}(z_t-\sigma_p)^2, \\
			\pi^c_t
			:&= c_t ( p_0 + p_1  S_t ) - \gamma c_t ( S_t + \delta ) - \frac{k_c}{2}v_t^2-\frac{\ell_c}{2}(y_t-\sigma_c)^2 ,
		\end{align*}
where $c_t ( p_0 + p_1  S_t )$ is the income from selling the quantity $c_t$ at the retail price $p_0 + p_1  S_t$, a linear function of the commodity price, with $p_0, p_1 > 0$ and $\gamma c_t (S_t + \delta )$ represents the sourcing cost of buying the quantity $\gamma c_t$ (which is used to obtain $c_t$ to be sold) at price $ S_t$ plus the transformation cost $\delta$, with $\gamma, \delta  > 0$. We assume $\gamma  > p_1$ to ensure the concavity of the objective functional of the consumer (i.e. the processor cannot charge increasing prices to final consumers without seeing the demand decreasing).  

The producer and the consumer exchange  a forward contract of $\lambda$ units of the commodity at a fixed amount of money $F \in \mathbb{R}$.  Both players aim at maximizing their respective objective functionals, which have two components: an expected profit term and a penalisation term modelling the player risk aversion (more comments below). 
In formulae, they are given by
\begin{align}	
J_p^{\lambda, F} (u, z;v,y)
	& :=\mathbb E \big[  P^p_T \big]   -  \eta_p  \int_0^T \mathbb{V} \left[ \lambda  S_t \right] dt, \qquad  \eta_p >0, \label{def_I_lambda} \\
J_c^{\lambda, F}(v,y; u, z)
	& :=\mathbb E \big[  P^c_T \big]    -  \eta_c  \int_0^T \mathbb{V} \left[ \lambda  S_t \right] dt, \qquad  \eta_c >0,\label{del_I_lambda_c}
\end{align}
where $\mathbb V$ stands for the variance and the process $P^p_T$ (resp. $P^c_T$) represents the cumulative profit over the time period $[0,T]$ of the producer (resp. the consumer), i.e.
\begin{align}\label{def_p^p_t}
		P_T^p : =  \int_0^T \pi_t^p dt + F  - \lambda  S_T,  \quad
		P_T^c :=  \int_0^T \pi_t^c dt - F + \lambda  S_T. 
\end{align}
The set of admissible strategies for the players is given by $\mathcal{A}^2:= \mathcal{A} \times \mathcal{A}$, where 
$\mathcal{A} = L_{\mathbb{F}}^2(\Omega \times [0,T], \mathbb{R}^2)$.\\

The way risk aversion is modelled and the choice of the derivative require two comments. First, a more standard way to take into account the players' risk aversion would consist in using utility functions. In our case and with an exponential utility function, where players can control the volatility of their production and consumption rates, this approach would lead to Monge-Amp\`ere PDEs, which are difficult to handle. For this reason, we turn to a different way to model risk aversion, which is reminiscent of what is done in mean-variance optimal dynamic portfolio choice (see \cite{ZL00} and more recently by \cite{IP19} and \cite{LLP20}). A similar approach was also previously used for distributed renewable energy development in \cite{ABP20}. Second, we observe that the variance penalisation term involves only the derivative and not the profit from production or transformation. As already stated in the introduction, this representation of risk aversion transforms players into speculators on the forward market. Indeed, players only care about the variance of  their financial position $\lambda S_t - F$, not about their production or consumption profits. This modeling is motivated by the desire to remain in a framework where tractable solutions can be exhibited. Its sole consequence would be to reverse the sign of the risk premium: producers wish to sell at a lower price than the expected spot price whereas speculators want to sell at a higher price. For the sake of simplicity, we have chosen to consider only a static hedging position with a simple forward contract in order to analyse the risk premium between the forward "agreement indifference price" and the expected price at maturity (see Section~\ref{sect_indiff_price} for a definition of the forward agreement indifference price).\medskip

To sum up, we deal with a two-player stochastic differential game of McKean-Vlasov linear-quadratic type. Hence, it is natural to look for \emph{Nash equilibria} according to the following definition.

\begin{definition}
	We call the couple 
	$\left((u^*,z^*)^\top,(v^*, y^*)^\top\right) \in \mathcal A \times \mathcal A$ a \emph{Nash equilibrium} if
	\begin{align}
		&  J_p^{\lambda, F}(u^*,z^*;v^*,y^*)\geq  J_p	^{\lambda, F}(u,z;v^*,y^*), \qquad \text{ for all } (u,z)^\top \in \mathcal{A},\\&
		  J_c^{\lambda, F}(v^*,y^*;u^*,z^*)\geq  J_c^{\lambda, F}(v,y;u^*,z^*), \qquad \text{ for all } (v,y)^\top \in \mathcal{A}.
	\end{align} 
\end{definition}

\subsection{Equilibrium forward agreement} \label{sect_indiff_price}
	For a Nash equilibrium $(v^*,y^*;u^*,z^*)$, we denote by
\[
	J^* _c (\lambda, F)
	=J_c^{\lambda, F}(v^*,y^*;u^*,z^*), \qquad
	J^* _p (\lambda, F)
	=J_p^{\lambda, F}(u^*,z^*;v^*,y^*),
\]
the corresponding equilibrium payoffs of consumer and producer, respectively. They depend on the number of units $\lambda$, on which the forward contract is written, and the respective forward price $F$.
Both players determine their prices using the \emph{indifference pricing approach}, namely the consumer computes $F_c ^{\lambda,*}$ as solution of $J^* _c (\lambda, F) = J^* _c (0,0)$ and analogously for the producer, leading to a price $F_p ^{\lambda,*}$ as a solution of $J^* _p (\lambda, F) = J^* _p (0,0)$.
By linearity of the payoffs with respect to $F$, we get
\[ 
	 J^* _c(\lambda,F)
	= J^* _c(\lambda,0) - F
	\quad \text{and} \quad
	 J^* _p(\lambda,F)
	= J^* _p(\lambda,0) + F,
\]
yielding
\[ 
	F^{\lambda,*} _c 
	= J^* _c(\lambda,0) - J^*_c (0,0), 
	\quad \text{and} \quad 
	F^{\lambda,*} _p 
	= J^* _p(0,0) - J^* _p (\lambda,0).
\]
Thus, $ F_c ^{\lambda,\ast}$ represents the maximum amount the consumer is willing to pay, while $F_p ^{\lambda,\ast}$ is the minimum amount the producer is willing to accept for selling a forward contract on $\lambda$ units of the underlying.
As a consequence, trading is possible if and only if 
\begin{equation}\label{ineq_K_lambda}
	F_p ^{\lambda,\ast} \leq F_c ^{\lambda,\ast}.
\end{equation}
We conclude this section with the definition of agreement indifference price. 

\begin{definition}
		Let $\lambda^*$ be the number of units of the underlying for which the two parties agree on the forward price, namely $	F_p ^{\lambda^\ast,\ast} = F_c ^{\lambda^\ast,\ast} $.
		We define the \emph{agreement indifference price}  as
		\begin{equation*}
			{F}_{{\lambda}^\ast}^\ast
			:=F_p ^{\lambda^\ast,\ast}
			=F_c ^{\lambda^\ast,\ast}.
		\end{equation*}
	\end{definition}
In Section \ref{sec:riskpremia}, we will provide some numerical illustrations on how the risk aversion parameters and the volatility control costs of the players might affect the quantity $\lambda^*$ as well as the corresponding agreement indifference price ${F}_{{\lambda}^\ast}^\ast$.

\section{Nash equilibrium}\label{sec:nash}
In this section we state and comment the main result of the paper. In particular we show that a Nash equilibrium exists and we characterize the corresponding strategies and payoffs in a semi-explicit way. Its proof will be given in full detail in the next section.
 
\subsection{Main result}\label{sect_main_res}

Let us start with some useful notation: for $t \in [0,T]$,
\begin{align}\label{Ks_Lambdas}
	&K_p(t)=-\frac{k_p}{2}\sqrt{\frac{2(\rho_p+\eta_p \lambda^2 \rho_p^2)}{k_p}}\tanh \left( \sqrt{\frac{2(\rho_p+\eta_p \lambda^2 \rho_p^2)}{k_p}} (T-t)\right),\\&
	\Lambda_p(t)= -\frac{k_p}{2}\sqrt{\frac{2\rho_p}{k_p}}\tanh \left( \sqrt{\frac{2 \rho_p }{k_p}} (T-t)\right),\\&
	K_c(t)=-\frac{k_c}{2}\sqrt{\frac{2(\gamma \rho_c(\gamma-p_1)+\eta_c  \lambda^2 \gamma^2\rho_c^2)}{k_c}}\tanh \left( \sqrt{\frac{2(\gamma \rho_c(\gamma-p_1)+\eta_c \lambda^2 \gamma^2\rho_c^2)}{k_c}} (T-t)\right),\\&
	\Lambda_c(t)= -\frac{k_c}{2}\sqrt{\frac{2\gamma\rho_c(\gamma-p_1)}{k_c}}\tanh \left( \sqrt{\frac{2 \gamma\rho_c(\gamma-p_1) }{k_c}} (T-t)\right),
\end{align}

\begin{align*}	
	\nonumber
	&\Xi=
	\begin{pmatrix}
		0 & -\rho_p\gamma \rho_c \eta_p \lambda^2-\frac{\gamma \rho_c}{2} \\
		-\rho_p\gamma \rho_c \eta_c \lambda^2 -\frac{\rho_p(\gamma-p_1)}{2}& 0 
	\end{pmatrix}, \quad
	\widehat\Xi=
	\begin{pmatrix}
		0 & -\frac{\gamma \rho_c}{2} \\
		-\frac{\rho_p(\gamma-p_1)}{2}& 0 
	\end{pmatrix}, \quad
	R =
	\begin{pmatrix}
		-\frac{2}{k_p} & 0 \\
		0 & -\frac{2}{k_c} 
	\end{pmatrix},
\end{align*}
\begin{align*}
	\Phi(t) =
	\begin{pmatrix}
		-\frac{2}{k_p}K_p(t) & 0 \\
		0 & -\frac{2}{k_c}K_c(t)  
	\end{pmatrix},\quad
	\widehat \Phi(t) =
	\begin{pmatrix}
		-\frac{2}{k_p}\Lambda_p(t) & 0 \\
		0 & -\frac{2}{k_c}\Lambda_c(t)  
	\end{pmatrix},\quad
	 \Psi= 
	\begin{pmatrix}
		-s_0/2 \\ 
		-\frac{p_0+p_1 s_0-\gamma(\delta+s_0)}{2}
	\end{pmatrix}.
\end{align*}
Furthermore, let us introduce the following system of ODEs  defined on $t \in [0,T]$:
\begin{equation}\label{Sys_Nash_simple_1}
	\begin{split}
		\left\{\begin{array}{ll}
			\pi'(t)=\Xi+\Phi(t) \pi(t)+\pi(t)\Phi(t)+\pi(t)R\pi(t),&
			\pi(T)=0,\\
			\widehat\pi'(t)=\widehat\Xi+\widehat\Phi(t) \widehat\pi(t)+	\widehat\pi(t)\widehat\Phi(t)+\widehat\pi(t)R\widehat\pi(t),&
			\widehat\pi(T)=0,
		\end{array}\right.
	\end{split}
\end{equation}
\begin{equation}\label{Sys_Nash_simple_2}
	\begin{split}
		\begin{array}{ll}
		dh(t)=\Big\{\big[\widehat\pi(t)R+\widehat\Phi(t)\big]h(t)+\Psi\Big\}dt,&
		h(T)=\frac{1}{2}\lambda(\rho_p,\gamma \rho_c)^\top,
		\end{array}
	\end{split}
\end{equation}
and let us denote by $T_{max}$ the right end of the interval where the system \eqref{Sys_Nash_simple_1} admits a unique solution according to Picard-Lindel\"{o}f Theorem (see, e.g., \cite{Coddington}, Theorem 2.3, which can be applied by standard time-inversion). 

\begin{theorem}\label{thm_Nash}
	Assume that the following conditions hold: 
	\begin{itemize}
		\item[\textbf{(A1)}] $T < T_{max}$,
		\item[\textbf{(A2)}] 
		$\ell_p-2(K_p(t)+\pi_{11}(t)) >0$ and $\ell_c-2(K_c(t)+\pi_{22}(t))>0$, for all $t \in [0,T]$.
	\end{itemize}	
	Then, 
	\begin{enumerate}
	\item there exists a Nash equilibrium  
	$((u^*,z^*)^{\top},(v^*,y^*)^\top) \in \mathcal{A}^2$ in the following feedback form
	\begin{align}
		&u^*_t = \frac{2}{k_p} \Big[(K_p(t)+\pi_{11}(t))(q_t-\bar q_t)+\pi_{12}(t)(c_t-\bar c_t)+(\Lambda_p(t)+\widehat \pi_{11}(t))\bar q_t+\widehat\pi_{12}(t)\bar c_t +h_1(t)\Big],\\ &\nonumber
		z^*(t) = \frac{\sigma_p \ell_p}{\ell_p-2(K_p(t)+\pi_{11}(t))},\\&\nonumber
		v^*_t = \frac{2}{k_c} \Big[(K_c(t)+\pi_{22}(t))(c_t-\bar c_t)+\pi_{21}(t)(q_t-\bar q_t)+(\Lambda_c(t)+\widehat \pi_{22}(t))\bar c_t+\widehat\pi_{21}(t)\bar q_t +h_2(t) \Big],\\&\nonumber
		y^*(t) = \frac{\sigma_c \ell_c}{\ell_c-2(K_c(t)+\pi_{22}(t))}.
	\end{align}
	\item The equilibrium payoffs satisfy  
	\begin{align}
	J_p^*(\lambda,F) & = \Lambda_p(0) q_0^2 + 2 \bar Y^p_0 q_0 + R_p(0) + F - \lambda s_0 - \frac12 \ell_p \sigma_p^2 T, \label{eq:Jp} \\
	J_c^*(\lambda,F) & =  \Lambda_c(0) c_0^2 + 2 \bar Y^c_0 c_0 + R_c(0)  - F + \lambda s_0 - \frac12 \ell_c \sigma_c^2 T, \label{eq:Jc}\\
	\bar Y^p_t & = \widehat \pi_{11} \bar q_t + \widehat \pi_{12} \bar c_t + h_1, \quad 
	\bar Y^c_t = \widehat \pi_{22} \bar c_t + \widehat \pi_{21} \bar q_t + h_2,
	\end{align}	
	where
	\begin{align}\label{eqs:R_expl}
		R_p(0)=R^{(\lambda)}_p(0)=& \int_0^T \left[ \frac{2}{k_p} \E [(Y^p_u)^2]   -\eta_p \lambda^2 \gamma^2\rho_c^2  \V [c_u] 
 +\frac{2\big(     \pi_{11}(u) z^\ast_u  + \frac{\ell_p \sigma_p}{2}\big )^2}{\ell_p-2K_p(u)}  \right] du - \lambda \gamma \rho_c \E [c_T], \\
		\nonumber R_c(0)=R^{(\lambda)}_c(0)= &   \int_0^T \left[   \frac{2}{k_c}  \E [(Y^c_u)^2]  
-\eta_c \lambda^2 \rho_p^2 \V [q_u]
 + \frac{2\big( \pi_{22}(u) y^\ast_u + \frac{\ell_c \sigma_c}{2}\big)^2}{\ell_c-2K_c(u)}   \right ]du - \lambda \rho_p \E [q_T].
	\end{align}
	
	\end{enumerate}
\end{theorem}

See Appendix~\ref{app:valfunc} for the details on the computations of the quantities involved in the definition of 
$R_p$ and $R_c$. 

\subsection{Comments}

\begin{enumerate}
\item Despite our model is close to the one presented in \cite{MP18}, it is not possible to directly exploit their results, since their hypotheses \textbf{(H2)(a)}  and \textbf{(H2)(d)} are not satisfied in our case. 
Therefore, in order to be self contained, we decided to prove a suitable verification theorem from scratch. 

\item We observe that the functions $\Lambda_{i}$, $i \in\{ p, c \}$, do not depend on $\lambda$. It is also the case for the functions $\widehat \pi_{ij}$. Furthermore, the functions $h_i$, $i=1,2$, are linear in $\lambda$ because they depend on it only by their terminal conditions. Besides, they are also nondecreasing in $\lambda$. Thus, the average production and consumption rates, $\bar q_t$ and $\bar c_t$, which satisfy
\begin{align*}
d\bar q_t &= \bar u^\ast_t dt =  \frac{2}{k_p} \Big[ \big( \Lambda_p(t) +\widehat \pi_{11}(t) \big) \bar q_t  + \widehat \pi_{12}(t)  \bar c_t +  h_1(t) \Big]dt , \\
d\bar c_t &= \bar v^\ast_t dt =  \frac{2}{k_c}  \Big[ \big( \Lambda_c (t) +\widehat \pi_{22} (t) \big) \bar c_t  + \widehat \pi_{21} (t)   \bar q_t +  h_2 (t)   \Big]dt, 
\end{align*}
are increasing in $\lambda$. As the terminal conditions of $h_{i}$, $i\in \{p,c\}$, depend only on  $\lambda$ and on the parameters $\rho_p$ and $\gamma \rho_c$, the resulting effect on the average spot price $\bar S_t = s_0 - \rho_p \bar q_t + \gamma \rho_c \bar c_t$ only depends on the relative market power of the producer and the consumer.  Thus, if $\gamma \rho_c < \rho_p$ (resp. $\rho_p < \gamma \rho_c$),  when the quantity of the commodity $\lambda$ of the producer increases, the average spot price decreases (resp. increases).

\item The functions $\Lambda_p$, $\Lambda_c$ and the $\widehat \pi_{ij}$ do not depend on the risk aversion parameters $\eta_p$ and $\eta_c$, therefore the average production and consumption rates do not depend on them either, as one could expect. Regarding the volatilities, while it is clear that $K_p$ and $K_c$ are nondecreasing in $\eta_p$ and $\eta_c$, respectively, it is not so obvious what to expect for $\pi_{11}$ and $\pi_{22}$, and thus to deduce the effect of risk-aversion on $z^\ast$ and $y^\ast$.  However, one can find numerically that the higher the risk aversions of the players, the lower the volatilities, even in the absence of forward agreement.  Nevertheless, it is possible to provide more insight on this issue when the producer has no market power, i.e. $\rho_p=0$, and the consumer does have some, i.e. $\gamma \rho_c>0$. In this case, the price process appears as exogenously driven for the producer and as a controlled variable for the consumer.  Hence $K_p = \Lambda_p = 0$ and $K_c <0$, $\Lambda_c<0$. Further, if $\rho_p=0$, then $\pi_{21} = 0$, leading to $\pi_{11} = 0$ due to $K_p = 0$, and it holds also that $\pi_{22}=0$  and $\widehat \pi_{11}= \widehat{\pi}_{21}= 0$.  Thus, $z^\ast = \sigma_p$  and the producer does not reduce her volatility. On the other hand, the production does covariate with consumption. Indeed, in Theorem~\ref{thm_Nash}, the Nash equilibrium consumer's strategies depend  on the state variables only via $c_t - \bar c_t$ and $\bar c_t$:
 \begin{align*}
&  u^\ast_t  = \frac{2}{k_p}\left\{\pi_{12}(t) (c_t - \bar c_t ) + \widehat \pi_{12} (t) \bar c_t + h_1 (t) \right\}, \quad z^\ast_t = \sigma_p, \\
& v^\ast_t = \frac{2}{k_c} \left\{ K_c (t)  (c_t - \bar c_t) + \Lambda_c (t)\bar c_t + h_2 (t) \right\}, 
\quad y^\ast_t = \frac{\sigma_c \ell_c}{\ell_c - 2 K_c(t)} < \sigma_c.
 \end{align*}
Finally, since $K_c(t)$ is nonincreasing in $\lambda$, the higher the exposure to  the financial risk coming from the forward contract, the more the consumer reduces his volatility, as the intuition predicts. 
\item Exploiting Theorem \ref{thm_Nash}-\emph{2.}, we can specify more precisely the nonlinear equations giving the forward agreement values $F^\ast_{{\lambda}^\ast}$ and $\lambda^\ast$. 
Indeed, it holds that (see equations \eqref{eq:Jp} and \eqref{eq:Jc})
\begin{align*}
&J_p^\ast(\lambda,F)  = \Lambda_p(0) q_0^2 + 2 \bar Y^{p (\lambda)}_0 q_0 + R_p^{(\lambda)}(0) + F - \lambda s_0 - \frac12 \ell_p \sigma_p^2 T, \\
&J_c^\ast(\lambda,F)  =  \Lambda_c(0) c_0^2 + 2 \bar Y^{c(\lambda)}_0 c_0 + R_c^{(\lambda)}(0)  - F + \lambda s_0 - \frac12 \ell_c \sigma_c^2 T,
\end{align*}
where the superscript $(\lambda)$ is used to emphasize the dependency on the number of options traded. We can isolate the parts $j_p(\lambda,F)$ and $j_c(\lambda,F)$ depending on $\lambda$ and $F$, defined as
\begin{align*}
&j_p(\lambda,F)  = 2 h_1^{(\lambda)}(0) q_0 + R^{(\lambda)}_p(0) + F - \lambda s_0, \\
& j_c(\lambda,F)  =   2 h_2^{(\lambda)}(0)  c_0 + R^{(\lambda)}_c(0)  - F + \lambda s_0.
\end{align*}
Thus, for a fixed $\lambda$ the indifference prices $F_p^{\lambda,\ast}$ and $F_c^{\lambda,\ast}$ are given by
\begin{align*}
& 2 h_1^{(0)}(0) q_0 + R^{(0)}_p(0)  = 2 h_1^{(\lambda)}(0) q_0 + R^{(\lambda)}_p(0) + F_p^{\lambda,\ast} - \lambda s_0, \\
& 2 h_2^{(0)}(0)  c_0 + R^{(0)}_c(0)   =   2 h_2^{(\lambda)}(0)  c_0 + R^{(\lambda)}_c(0)  - F_c^{\lambda,\ast} + \lambda s_0.
\end{align*}
Thus, if it exists, the equilibrium price should be given by $F_p^{\lambda^\ast,\ast}$ $=$ $F_c^{\lambda^\ast,\ast}={F}_{{\lambda}^\ast}^\ast$, i.e.,
\begin{align*}
& 2 (h_1^{(0)}(0) -h_1^{(\lambda^\ast)}(0)) q_0 + R^{(0)}_p(0) - R^{(\lambda^\ast)}_p(0) =
 2 (h_2^{(\lambda^\ast)}(0) - h_2^{(0)}(0) )  c_0 + R^{(\lambda^\ast)}_c(0) - R^{(0)}_c(0), 
\end{align*}
or, equivalently, 
\begin{align}
& 2 h_1^{(0)}(0) q_0 + 2 h_2^{(0)}(0) c_0 + R^{(0)}_c(0) + R^{(0)}_p(0)  =
 2 h_1^{(\lambda^\ast)}(0) q_0 + 2 h_2^{(\lambda^\ast)}(0) c_0 + R^{(\lambda^\ast)}_c(0) + R^{(\lambda^\ast)}_p(0),
\end{align}
with $R^{(\lambda^\ast)}_p(0)$ and $R^{(\lambda^\ast)}_c(0)$ defined in  Equation \eqref{eqs:R_expl} and $h^{(\lambda^\ast)}$ in Equation \eqref{Sys_Nash_simple_2}.
\end{enumerate}

The last remark speeds up considerably the computations for the plots that appear in the Section \ref{sec:numeric}. Indeed, all the quantities that we need to compute can be obtained by solving numerically the ODEs presented in Appendix \ref{app:valfunc}.

\section{Proof of Theorem \ref{thm_Nash}}\label{sec:proof}

\subsection{The solution approach}
We prove Theorem \ref{thm_Nash} following a methodology based on a combination of a suitable \emph{Verification Theorem} and of the \emph{weak Martingale Optimality Principle}. As already stressed in the first comment below Theorem \ref{thm_Nash}, despite our model is very close to the class of games studied in \cite{MP18}, their results cannot be applied directly here, therefore we had to adapt the methodology to our framework. We proceed through the following steps:
\begin{itemize}
	\item[1)]  we compute the best response maps of both players;

	\item[2)] we check that the system coming from the best response computations has a unique solution;

	\item[3)] we get a Nash equilibrium as a fixed point of the best response map;

	\item[4)] we verify that there exists a unique solution to the system characterizing the fixed point found in step 3).
\end{itemize}

\subsection{Preliminary reformulation of the problem}
For convenience, we introduce the following vector notation for the players' strategies:
\begin{align*}
	\alpha=\left((\alpha^p)^\top,(\alpha^c)^\top \right)^\top \in \mathcal{A}^2, \quad 
	\alpha^p:= 
	\begin{pmatrix}
		u  \\
		z
	\end{pmatrix}
	=\left\{
	\begin{pmatrix}
		u_t  \\
		z_t
	\end{pmatrix}\right\}_{t \in [0,T]}
	\text{ and }\quad
	\alpha^c:= 
	\begin{pmatrix}
		v  \\
		y
	\end{pmatrix}
	= \left\{
	\begin{pmatrix}
		v_t  \\
		y_t
	\end{pmatrix}\right\}_{t \in [0,T]},
\end{align*}
so that the dynamics of the state variables can be rewritten as 
\begin{align}\label{eq_state_game}
	&dq_t=e_1^{\top} \alpha^p_t dt + e_2^\top \alpha^p_t dW_t,\\&
	dc_t=e_1^\top \alpha^c_t dt + e_2^\top \alpha^c_t dB_t, \qquad t \in [0,T],
\end{align}
with
$
	e_1^\top = (1, 0)
	\text{ and }
	e_2^\top = (0, 1).
$
\\

The following identity is exploited to get a suitable reformulation of our problem: using the dynamics of $S_t$ and applying Fubini's theorem, it is easy to see that
\begin{align}\label{variance_id}
	\int_0^T \mathbb V \left[ S_t \right] dt 
	= \mathbb E \left[ \int_0^T \left\{ \rho_p^2 (q_t - \mathbb E [q_t])^2 + \gamma^2\rho_c^2 (c_t - \mathbb E [c_t])^2 -2 \rho_p \gamma \rho_c (q_t - \mathbb E [q_t])(c_t - \mathbb E [c_t]) \right\} dt\right].
\end{align}
Rearranging the terms in the expressions of the producer objective functional, we obtain
\begin{displaymath}
	\begin{split}
		J_p^{\lambda, F}(u, z; v, y) &
		= \widetilde{J}_p^{\lambda}(u,z;v,y) +F -\lambda s_0 -\frac{\ell_p \sigma_p^2 T}{2},
	\end{split}
\end{displaymath}
where
\begin{equation}\label{eq_def_widetilde_I}
	\begin{split}
		\widetilde{J}_p^{\lambda}(u,z;v,y) 
		:= & \mathbb{E} \Bigg[ \int_0^T\bigg(  -(\rho_p +\eta_p \lambda^2 \rho_p^2)(q_t-\mathbb E [q_t])^2 - \rho_p \mathbb E [q_t]^2+[s_0+\gamma \rho_c c_t\\&
		+2\rho_p \gamma \rho_c \eta_p\lambda^2 (c_t- \mathbb E [c_t])] q_t - \frac{k_p}{2}u_t^2-\frac{\ell_p}{2}z_t^2 + \ell_p \sigma_p z_t -\eta_p\lambda^2 \gamma^2\rho_c^2 (c_t -\mathbb E [c_t])^2 \bigg)  dt \\&
		+\lambda \rho_p q_T -\lambda \gamma \rho_c c_T \Bigg].
	\end{split}
\end{equation}
Then, neglecting the constant terms, we can study without loss of generality the equivalent formulation in which the producer aims at maximizing $\widetilde{J}_p^{\lambda}(u,z;v,y) $.

\begin{remark}\label{rem_fix_con_fix_stat}
	Fixing a strategy $\alpha^p$ for the producer (resp. $\alpha^c$ for  the consumer) is equivalent, from the perspective of the competitor, to fixing the corresponding state $q^{\alpha^p}$ (resp. $c^{\alpha^c}$). Thus, with some abuse of notation we will write simply $q$ (resp. $c$) when the strategy used is clear from the context. Moreover, to ease the notation, we will also omit the dependence on $\bar c$ and $\bar q$.
\end{remark}

For a given consumption process $\{c_t\}_{t \in [0,T]}$, we write

\begin{equation}\label{def_widehat_I_c}
	\begin{split}
		\widetilde J_p^\lambda& (\alpha^p;\alpha^c)
		=\widetilde J_p^\lambda (\alpha^p;c)
		:=\mathbb E\left[ \int_0^T f_p(t,q_t,\mathbb E [q_t], \alpha^p_t, \mathbb E[\alpha^p_t];c) dt + g_p(q_T,\mathbb E [q_T];c)  \right], \text{ with }\\&
		f_p(t,q,\bar{q}, a_p, \bar{a}_p;c)=   Q_p (q-\bar{q})^2 + (Q_p+\widetilde Q_p)\bar{q}^2 +2 M^{p}(c)_t q+a_p^\top N_p a_p +2H_p^\top a_p +T^{p}(c)_t, \\&
		g_p(q,\bar{q};c) =   2L_pq + \widetilde T^{p}(c) ,
	\end{split}
\end{equation}
where
\begin{equation}\label{eq_matrices_p}
	\begin{split}
		&Q_p := -\rho_p - \eta_p \lambda^2 \rho_p^2, \quad 
		\widetilde Q_p:=  \eta_p \lambda^2 \rho_p^2,
		\quad
		M^{p}(c)_t := \frac{s_0}{2} + \frac{\gamma \rho_c}{2}c_t + \rho_p \gamma \rho_c \eta_p \lambda^2 (c_t - \mathbb E [c_t]),
		\\&
		N_p :=
		\begin{pmatrix}
			-\frac{k_p}{2} & 0\\
			0 & -\frac{\ell_p}{2}
		\end{pmatrix},
		\quad
		H_p := 
		\begin{pmatrix}
			0   \\
			\frac{\sigma_p \ell_p}{2}
		\end{pmatrix}, \quad
		T^{p}(c)_t:= -\eta_p \lambda^2 \gamma^2\rho_c^2 (c_t-\mathbb E[c_t])^2, 
		\quad
		L_p := \frac{\rho_p \lambda}{2},\\& 
		\text{ and }\quad 
		\widetilde T^{p}(c):=-{\lambda \gamma \rho_c}c_T.
	\end{split}
\end{equation}
Now, let us turn to the objective functional of the consumer. From \eqref{del_I_lambda_c} and \eqref{variance_id}, we have
\begin{displaymath}
	\begin{split}
		{J}_c^{\lambda, F} & (v, y; u, z) 
		= \widetilde{J}_c^{\lambda}(v,y;u,z)-F + \lambda s_0  -\frac{\ell_c \sigma_c^2T}{2},
	\end{split}
\end{displaymath}
where 
\begin{equation}
	\begin{split}
		\widetilde{J}_c^{\lambda}(v,y; u,z)
		:=  \mathbb{E}\Bigg[ & \int_0^T\bigg(  -[\gamma \rho_c(\gamma-p_1)+\eta_c \lambda^2 \gamma^2\rho_c^2](c_t-\mathbb E [c_t])^2 - \gamma \rho_c(\gamma-p_1)\mathbb E[c_t]^2  \\& 
		+[ (p_0+s_0 p_1 -\gamma \delta-\gamma s_0)  + \rho_p(\gamma-p_1) q_t+2 \rho_p \gamma \rho_c \eta_c \lambda^2 (q_t-\mathbb E [q_t]) ] c_t \\&
		 - \frac{k_c}{2}v_t^2-\frac{\ell_c}{2}y_t^2 + \sigma_c \ell_c y_t  - \eta_c \lambda^2 \rho_p^2 (q_t - \mathbb E [q_t])^2 \bigg)  dt 
		-\lambda \rho_p q_T +\lambda \gamma \rho_c c_T\Bigg]. 
	\end{split}
\end{equation}
Analogously as above, let $\{q_t\}_{t \in [0,T]}$ be a given production rate. We can write 
\begin{equation*}
	\begin{split}
		\widetilde J_c^\lambda&(\alpha^c; \alpha^p)
		=\widetilde J_c^\lambda(\alpha^c; q)
		:=\mathbb E\left[ \int_0^Tf_c\big(t,c_t,\mathbb E [c_t], \alpha^c_t, \mathbb E[\alpha^c_t];q\big) dt + g_c(c_T,\mathbb E [c_T];q)  \right], \text{ with }\\&
		f_c(t,c,\bar{c}, a_c, \bar{a}_c;q)=   Q_c (c-\bar{c})^2 + (Q_c+\widetilde Q_c)\bar{c}^2 +2 M^{c}(q)_t c+a_c^\top N_c a_c +2H_c^\top a_c + T^{c}(q)_t, \\&
		g_c(c,\bar{c};q) =   2L_cc + \widetilde T^{c}(q) ,
	\end{split}
\end{equation*}
and 
\begin{equation}\label{eq_matrices_c}
	\begin{split}
		&Q_c :=-\gamma \rho_c(\gamma-p_1)-\eta_c \lambda^2\gamma^2 \rho_c^2, \quad 
		\widetilde Q_c:=  \eta_c \lambda^2 \gamma^2 \rho_c^2,
		\quad
		N_c :=
		\begin{pmatrix}
			-\frac{k_c}{2} & 0\\
			0 & -\frac{\ell_c}{2}
		\end{pmatrix},\\&
		M^{c}(q)_t := \frac{p_0+p_1 s_0-\gamma(s_0+\delta)}{2}+		\frac{\rho_p(\gamma-p_1)}{2}q_t+\rho_p \gamma \rho_c \eta_c \lambda^2 (q_t - \mathbb E [q_t]), \\&
		H_c := 
		\begin{pmatrix}
			0   \\
			\frac{\sigma_c \ell_c}{2}
		\end{pmatrix}, \quad
		T^{c}(q)_t:= -\eta_c \lambda^2 \rho_p^2 (q_t-\mathbb E[q_t])^2, 
		\quad
		L_c := \frac{ \lambda \gamma \rho_c}{2}, 
		\quad \text{ and }\quad 
		\widetilde T^{c}(q):=-{\lambda \rho_p}q_T.
	\end{split}
\end{equation}
Finally, we set
\begin{align}
	&V_p^\lambda(\alpha^c)
	:=\sup_{\alpha^p \in \mathcal{A}}\widetilde J_p^\lambda(\alpha^p; \alpha^c),
	\qquad \alpha^c \in \mathcal{A},\\&
	V_c^\lambda(\alpha^p)
	:=\sup_{\alpha^c \in \mathcal{A}}\widetilde  J_c^\lambda(\alpha^c; \alpha^p),
	\qquad \alpha^p \in \mathcal{A}.
\end{align}

\subsection{First step: computation of the best response maps}

The first step is focused on the computation of the best response map of each player. This is done by exploiting the following version of the \emph{Verification Theorem}: 

\begin{theorem}[Verification Theorem]\label{Ver_thm}
	Fix a couple of strategies $\beta^p, \beta^c \in \mathcal{A}$ for the producer and the consumer, respectively. Let $\mathcal{W}^{p,\alpha^p}_t$ and $\mathcal{W}^{c,\alpha^c}_t$  be defined as 
	\begin{align}
		&\mathcal{W}^{p,\alpha^p}_t= w_t^p(q_t^{\alpha^p}, \mathbb{E}[q_t^{\alpha^p}]) ,
		\quad 
		\mathcal{W}^{c,\alpha^c}_t= w_t^c(c_t^{\alpha^c}, \mathbb{E}[c_t^{\alpha^c}]), 
		\quad t \in [0,T], \quad \alpha^p, \alpha^c \in \mathcal{A},
	\end{align}
	where the $\mathbb{F}$-adapted random fields $\{w_t^p(q,\bar q), t \in [0,T], q, \bar q \in \mathbb{R}\}$ and $\{w_t^c(c,\bar c), t \in [0,T], c, \bar c \in \mathbb{R}\}$  satisfy the following growth conditions: $\text{ for all } t \in [0,T], \text{ for all } x,\bar{x} \in \mathbb{R},$
	\begin{align}\label{growth_cond}
		&|w_t^p(x,\bar x)| \leq C_p(\nu_t^p+|x|^2+|\bar{x}|^2),
		\qquad
		|w_t^c(x,\bar x)| \leq C_c(\nu_t^c+|x|^2+|\bar{x}|^2),
	\end{align}
	for some constants $C_p, C_c > 0$ and for some non-negative processes $\nu^p$ and $ \nu^c$ such that
	\[
		\sup_{t \in [0,T]}\mathbb{E}\left[\nu_t^p+\nu_t^c\right]< \infty.
	\]
	Furthermore, we assume that the following conditions are fulfilled:
	\begin{itemize}
		\item[i)] $\mathbb{E}[w_T^p(q_T^{\alpha^p},\bar{q_T}^{\alpha^p})]= \mathbb{E}[g_p(q_T^{\alpha^p},\bar{q_T}^{\alpha^p};c^{\beta^c})]$ and $\mathbb{E}[w_T^c(c_T^{\alpha^c},\bar{c_T}^{\alpha^c})]= \mathbb{E}[g_c(c_T^{\alpha^c},\bar{c_T}^{\alpha^c};q^{\beta^p})]$, for any $\alpha^p, \alpha^c \in \mathcal{A}$.
		
		\item[ii)] The application $[0,T] \ni t \mapsto \mathbb{E}[\mathcal{S}^{p,{\alpha^p}}_t] \Big(\text{resp. } \mathbb{E}[\mathcal{S}^{c,{\alpha^p}}_t]\Big)$ is well-defined and nonincreasing, for any $ \alpha^p \in \mathcal{A}$ (resp. for any $ \alpha^c \in \mathcal{A}$), where:
		\begin{equation}\label{def_mathcal_S}
			\begin{split}
				&\mathcal{S}^{p,{\alpha^p}}_t = \mathcal{W}^{p,\alpha^p}_t + \int_0^t f_p(s,q_s^{\alpha^p},\bar q_s^{\alpha^p}, \alpha^p_s, \bar \alpha^p_s; c^{\beta^c} )ds,\\&
				\mathcal{S}^{c,{\alpha^p}}_t = \mathcal{W}^{c,\alpha^c}_t + \int_0^t f_c(s,c_s^{\alpha^c},\bar c_s^{\alpha^c}, \alpha^c_s, \bar \alpha^c_s; q^{\beta^p} )ds.
			\end{split}
		\end{equation}
		
		\item[iii)] For some $\alpha^{p,\star} \in \mathcal{A}$ and $\alpha^{c,\star} \in \mathcal{A}$, the application $[0,T] \ni t \mapsto \mathbb{E}[\mathcal{S}^{p,\alpha^{p,\star}}_t] \Big(\text{resp. } \mathbb{E}[\mathcal{S}^{c,\alpha^{c,\star}}_t]\Big)$ is constant.	
	\end{itemize}
	Then, the control $\alpha^\star=(\alpha^{p,\star},\alpha^{c,\star})$ is the \emph{best response} to the control $(\beta^p,\beta^c)$ meaning that 
	\begin{equation}
	\alpha^{p,\star}= \mathbf{B}_p(\beta^c):= \argmax_{\alpha^p \in \mathcal{A}}\widetilde J_p^\lambda(\alpha^p;\beta^c),
	\quad
	\alpha^{c,\star}= \mathbf{B}_c(\beta^p):= \argmax_{\alpha^c \in \mathcal{A}}\widetilde J_c^\lambda(\alpha^c;\beta^p),
	\end{equation}		
	
	 and 
	\[
		\widetilde{J}_p^\lambda(\alpha^{p,\star};c^{\beta^c})
		=V^\lambda_p(\beta^c)
		= \mathbb{E}[\mathcal{W}_0^{P, \alpha^{p}}] \quad and 
		\quad
		\widetilde{J}_c^\lambda(\alpha^{c,\star};q^{\beta^p})
		=V^\lambda_c(\beta^p)
		= \mathbb{E}[\mathcal{W}_0^{C, \alpha^{c}}].
	\]
	Finally, if $\widetilde \alpha= (\widetilde \alpha^p, \widetilde \alpha^c)$ is another best response to the control $(\beta^p,\beta^c)$, then condition iii) holds also for $\widetilde \alpha^p$ and $\widetilde \alpha^c.$
\end{theorem}

We define the best response map $\textbf{B}: \mathcal{A}^2 \to \mathcal{A}^2$ as $\textbf{B}:=(\textbf{B}_p, \textbf{B}_c)$. 
The Nash equilibrium we find will be a fixed point of this map.
\\

Once we have fixed the strategies $\beta^p $ and $\beta^c$ in $\mathcal{A}$, the first step can be divided into four sub-steps:
\begin{itemize}
	\item[\textbf{1.1}] Since the players objective functionals are quadratic, we propose a suitable candidate $(\mathcal{W}^{p,\alpha^p}_t,\mathcal{W}^{c,\alpha^c}_t)$ in feedback form.

	\item[\textbf{1.2}]Applying It\^o's formula, we compute $\frac {d}{dt} \mathbb E [\mathcal{S}^{p,\alpha^p}_t]$ and $\frac {d}{dt} \mathbb E [\mathcal S^{c,\alpha^c}_t]$  corresponding to the candidate $(\mathcal{W}^{p,\alpha^p}_t,\mathcal{W}^{c,\alpha^c}_t)$.

	\item[\textbf{1.3}] We postulate that the conditions of Theorem \ref{Ver_thm} are satisfied and get a system of backward SDEs involving the coefficients of the candidate $(\mathcal{W}^{p,\alpha^p}_t,\mathcal{W}^{c,\alpha^c}_t)$.

	\item[\textbf{1.4}] We compute each player's best response by looking for strategies cancelling the expectation of the drifts of the processes $\mathcal{S}^{p,\alpha^p}_t$ and $\mathcal{S}^{c,\alpha^c}_t$.
 
\end{itemize}

\paragraph{\textbf{Sub-step 1.1}} 

Given the quadratic nature of our objective functional, it seems natural to look for a family of processes $(\mathcal{W}^{p,\alpha^p}_t,\mathcal {W}^{c,\alpha^c}_t)_{t \in [0,T]}$ of the following form: $\mathcal{W}^{p,\alpha^p}_t= w^p_t(q_t^{\alpha^p}, \mathbb 	E[q_t^{\alpha^p}])$ and $\mathcal{W}^{c,\alpha^c}_t= w^c_t(c_t^{\alpha^c}, \mathbb 	E[c_t^{\alpha^c}])$, for some parametric adapted random field 
 $\{w_t^i(x, \bar x), t \in [0,T], x, \bar x \in \mathbb R\}, i \in \{p,c\},$ such that 
\[
	w_t^i(x, \bar x)=K_i(t)(x-\bar x)^2+\Lambda_i(t)\bar x^2+2Y^i_t x+ R_i(t),
\] 
with $(K_i,\Lambda_i,Y^i,R_i) \in L^{\infty}([0,T], \mathbb R_-)^2\times S^2_{\mathbb F}(\Omega \times [0,T], \mathbb R)\times L^{\infty}([0,T], \mathbb R)$, $ i \in \{p,c\}$, solving the systems of ODEs and SDEs:
\begin{equation}
	\begin{split}
		&\left\{\begin{array}{lll}
			dK_p(t)&=K_p'(t)dt, &
			K_p(T)=0,
			\\
			d\Lambda_p(t)&=\Lambda_p'(t)dt,&
			\Lambda_p(T)=0, 
			\\ 
			dY^p_t&=Y^p_t{}' dt+Z^{p,B}_tdB_t+Z^{p,W}_t dW_t,&
			Y^p_T=\frac{\lambda \rho_p}{2},
			\\
			dR_p(t)&=R_p'(t)dt, &
			R_p(T)=-\lambda \gamma \rho_c \mathbb E [c_T],
		\end{array}\right.
	\end{split}
\end{equation}

\begin{equation}
	\begin{split}
		&\left\{\begin{array}{lll}
			dK_c(t)&={K}_c'(t)dt, &
			K_c(T)=0,
			\\
			d\Lambda_c(t)&={\Lambda}_c'(t)dt,&
			\Lambda_c(T)=0, 
			\\ 
			dY^c_t&={Y^c_t}' dt+Z^{c,B}_tdB_t+Z^{c,W}_tdW_t&
			Y^c_t=\frac{\lambda \gamma \rho_c}{2},
			\\
			dR_c(t)&={R}_c'(t)dt, &
			R_c(T)= -\lambda \rho_p \mathbb E [q_T],
		\end{array}\right.
	\end{split}
\end{equation}
for some deterministic processes $K_i', \Lambda_i', R_i'$ and for some $\mathbb{F}$-adapted processes ${Y^i}', Z^{i,W}, Z^{i,B}$, $i \in \{p, c\}$.

\paragraph{\textbf{Sub-step 1.2}} \label{Sub-step_1.2}
For the sake of simplicity, from now on, we  explicitly develop  only the producer case. The consumer problem can be studied in same way. 
Let $t \in [0,T]$ and $\alpha^p \in \mathcal A.$ 
As in \eqref{def_mathcal_S} in Theorem \ref{Ver_thm} (Verification Theorem), we set
\[
	\mathcal{S}^{p,\alpha^p}_t= w_t^p(q_t^{\alpha^p}, \mathbb E[ q_t^{\alpha^p}])+   \int_0^t 	f_p(u,q_u^{\alpha^p}, \mathbb E[{q_u^{\alpha^p}}], 	\alpha_u^p, \mathbb E [{\alpha_u^p}]; c^{\beta^c}) du. 
\]
 In the following, we write simply $c$ instead of $c^{\beta^c}$ (resp. $q$ instead of $q^{\beta^p}$), when the strategies are clear from the context (see Remark \ref{rem_fix_con_fix_stat}).
After some computations (see Appendix \ref{Comp_for_game} for details), we obtain
\begin{equation}\label{eq:E(Sp)}
	\begin{split}
		\frac {d}{dt}\mathbb E[\mathcal{S}^{p,\alpha^p}_t]&
		=\mathbb E \Big[  
		(  {K}_p'(t) + Q_p)(q_t- \mathbb E[q_t])^2 
		+(  {\Lambda}_p'(t) +Q_p+ \widetilde Q_p)			\mathbb E[q_t]^2
		+2(  {Y^p_t}' + M^{p}(c)_t)q_t\\&
		\qquad + {R}_p'(t) + T^{p}(c)_t+ \chi^p_t(\alpha^p_t)
		\Big],
	\end{split}
\end{equation} 
where, for any $ t \in [0,T],$ we have set
\begin{equation}\label{S_p_def}
	\left \{ \begin{array}{l}
		\chi^p_t(\alpha^p_t):= (\alpha^p_t)^\top S_p(t)\alpha^p_t +2[U_p(t)(q_t-\mathbb E[q_t])+V_p(t)q_t+\xi^p_t+\bar{\xi}^p_t+O_p(t)]^\top	\alpha^p(t) \\
		S_p(t):=N_p+e_2K_p(t)e_2^\top\\
		U_p(t):=K_p(t)e_1\\
		V_p(t):=\Lambda_p(t)e_1\\
		O_p(t):=H_p+ e_1\mathbb{E}[Y^p_t]+e_2\mathbb{E}[Z^{p,W}_t]\\ 
		\xi^p_t:= H_p+ e_1 Y^p_t+e_2 Z^{p,W}_t\\
		\bar{\xi}^p_t:= H_p+ e_1\mathbb{E}[Y^p_t]+e_2\mathbb{E}[Z^{p,W}_t],
	\end{array}\right.
\end{equation}
where $Q_p, \widetilde{Q}_p, M^{p}(c), N_p, H_p$ and $T^{p}(c)$ are defined in Equation \eqref{eq_matrices_p}.
\\

\paragraph{\textbf{Sub-step 1.3}} 
Now, we find conditions granting that assumptions i), ii) and iii) of Theorem \ref{Ver_thm}, involving $\mathcal S^{p,\alpha^p}$, hold. 
Suppose that the matrix $S_p(t)$ is negative definite and thus invertible.
We check this later, verifying that $K_p(t)  \leq 0, \text{ for all } t \in [0,T]$ (see Remark \ref{K_and_Lambda}).  
We complete the squares and rewrite the equation \eqref{eq:E(Sp)} as
\[
	\begin{split}
		\frac {d}{dt}\mathbb E[\mathcal{S}^{p,\alpha^p}_t]&
		=\mathbb E \Big[ 
		\big(  {K}_p'(t) + Q_p- U_p(t)^\top S_p(t)^{-1}	U_p(t)\big)(q_t- \mathbb E[q_t])^2 \\&
		\qquad + \big(  {\Lambda}_p'(t) +Q_p+ \widetilde Q_p- V_p(t)^\top S_p(t)^{-1}V_p(t)\big) \mathbb E[q_t]^2\\&
		\qquad +2 \big[  {Y^p_t}' + M^{p}(c)_t- U_p(t)^\top S_p(t)^{-1}( \xi^p_t-\bar{ \xi}^p_t)-  V_p(t)^\top S_p(t)^{-1} O_p(t)\big]q_t\\&
		\qquad + {R}_p'(t) + T^{p}(c)_t - (\xi^p_t - \bar{\xi}^p_t)^\top S_p(t)^{-1}(\xi^p_t - \bar{\xi}^p_t) -  O_p(t)^\top 	S_p(t)^{-1} O_p(t)\\&
		\qquad + (\alpha^p_t-\eta^p_t)^\top S_p(t)^{-1} (\alpha^p_t-	\eta^p_t)
		\Big],
	\end{split}
\] 
where, for all $ t \in [0,T],$ we have defined
\begin{align}
	\eta^p_t:=-S_p(t)^{-1}\left[U_p(t)(q_t- \mathbb E[q_t])+V_p(t)\mathbb E[q_t]+(\xi^p_t-\bar{\xi}^p_t)+ O_p(t)\right].
\end{align}
Choosing  processes $K_p, \Lambda_p, Y^p$ and $R_p$, whose existence is shown in the next sub-step, that solve  the following system of BSDEs
\begin{equation}\label{Sys_p}
	\begin{split}
		&\left\{\begin{array}{l}
		\begin{array}{ll}
  			{K}_p'(t) + Q_p- U_p(t)^\top S_p(t)^{-1}U_p(t)=0,&
			K_p(T)=0,
		\end{array}
		\\
		\begin{array}{ll}
			 {\Lambda}_p'(t) +Q_p+ \widetilde Q_p- V_p(t)^\top S_p(t)^{-1}V_p(t) =0, &
			\Lambda_p(T)=0,
		\end{array}
		\\
		\begin{array}{ll}
			dY^p_t=& \left[  - M^{p}(c)_t+ U_p(t)^\top S_p(t)^{-1}(\xi^p_t - \bar{\xi}^p_t)+ V_p(t)^\top S_p(t)^{-1} O_p(t)\right]dt \\
			&+ Z^{p,B}_tdB_t + Z^{p,W}_t dW_t,\\
			Y^p_T=L_p,& 
		\end{array}
		\\
		\begin{array}{l}
 			{R}_p'(t) + \mathbb E [ T^{p}(c)_t- (\xi^p_t - \bar{\xi}^p_t)^\top S_p(t)^{-1}(\xi^p_t - \bar{\xi}^p_t) -  O_p(t)^\top S_p(t)^{-1} O_p(t) ] =0, \\
				R_p(T) = \mathbb{E}[\widetilde T^{p}(c)],
			\end{array}
		\end{array}\right.
	\end{split}
\end{equation}
we obtain
\begin{equation}\label{opt_cont_cond}
	\begin{split}
		\frac {d}{dt}\mathbb E[\mathcal{S}^{p,\alpha^p}_t]&
		=\mathbb E \Big[  (\alpha^p_t-\eta^p_t)^\top S_p(t)^{-1} (\alpha^p_t-\eta^p_t)
		\Big],
	\end{split}
\end{equation}
which is clearly non-positive for all $t \in [0,T]$,  since $S_p(t)$ (defined in Equation \eqref{S_p_def}) is negative definite for all $ t \in [0,T]$.

\begin{remark}\label{Dep_on_state}
	We stress the fact that the processes $Y^p, Z^{p,W},Z^{p,B}$ and $R_p$ depend only on the strategy of the consumer through the state process $\{c_t\}_{t \in [0,T]}$, with $c_t=c^{\beta^c}_t, t \in [0,T]$, which is controlled only by $\beta^c.$ Thus, the feedback best response control are functions of different state variables, namely the best response for the producer is feedback in $q$ and its expectation, whereas the best response for the consumer is feedback in $c$ and its expectation.
\end{remark}

\paragraph{\textbf{Sub-step 1.4}} 

Now we combine the results in the previous steps in order to get the best response maps.

\begin{proposition}\label{Prop_best_resp}
	The best response maps
	are given by
	\begin{align}\label{Best_r}
		&\mathbf{B}_p(\beta^c)_t
		= -(N_p+e_2 K_p(t)e_2^\top)^{-1}[e_1 K_p(t)(q_t-\mathbb E [q_t])+e_1\Lambda_p(t) \mathbb E [q_t]+e_1 Y^p_t+e_2 Z^{p,W}_t + H_p ],\\&\nonumber
		\mathbf{B}_c(\beta^p)_t
		= -(N_c+e_2 K_c(t)e_2^\top)^{-1}[e_1 K_c(t)(c_t-\mathbb E [c_t])+e_1\Lambda_c(t) \mathbb E [c_t]+e_1 Y^c_t+e_2 Z^{c,B}_t +H_c],
	\end{align}
	where the processes $(K_p,\Lambda_p,Y^p,R_p) $ and $(K_c,\Lambda_c,Y^c,R_c)$ above solve the following systems of backward ODEs and SDEs, given $c_t=c_t^{\beta^c}$ (respectively, given $q_t=q_t^{\beta^p}$), $t \in[0,T]$:
	\begin{equation}\label{Sys_p}
	\begin{split}
		&\left\{\begin{array}{l}
			\begin{array}{ll}
				K_p'(t) = -\frac{2}{k_p}K_p(t)^2+\rho_p + \eta_p \lambda^2 \rho_p^2,&
				K_p(T)=0,
			\end{array}
			\\
			\begin{array}{ll}
				\Lambda_p'(t) = -\frac{2}{k_p} 	\Lambda_p(t)^2 +\rho_p, &
				\Lambda_p(T)=0,
			\end{array}
			\\
			\begin{array}{ll}
				dY^p_t&=-\left\{\frac{s_0}{2}+\frac {\gamma \rho_c} {2} c_t +\rho_p \gamma \rho_c \eta_p \lambda^2 (c_t-\mathbb E [c_t])+\frac{2}{k_p}\Big[K_p(t)\left(Y^p_t  - \mathbb E [Y^p_t]\right)+\Lambda_p(t)\mathbb E[Y^p_t]\Big] \right\}dt\\
				&\quad +Z^{p,B}_tdB_t+Z^{p,W}_tdW_t,\\
				Y^p_T&=\frac{\lambda \rho_p}{2},
			\end{array}
			\\
			\begin{array}{l}
				R_p'(t)  = \eta_p \lambda^2 \gamma^2\rho_c^2 \mathbb V [c_t] -\frac{2}{k_p}(\mathbb V [Y^p_t]+\mathbb E [Y^p_t]^2)-\frac{2}{\ell_p-2K_p(t)}(\mathbb V[Z^{p,W}_t]+(\mathbb E [Z^{p,W}_t]+\frac{\ell_p 	\sigma_p}{2})^2), \\
				R_p(T)= - \lambda \gamma \rho_c \mathbb E [c_T],
			\end{array}
		\end{array}\right.
	\end{split}
	\end{equation}
and
	\begin{equation}\label{Sys_c}
	\begin{split}
		&\left\{\begin{array}{l}
			\begin{array}{ll}
				K_c'(t) = -\frac{2}{k_c}K_c(t)^2+\gamma \rho_c(\gamma-p_1)+\eta_c \lambda^2 \gamma^2 \rho_c^2=0,&
				K_c(T)=0,
			\end{array}
			\\
			\begin{array}{ll}
				\Lambda_c'(t) = -\frac{2}{k_c} 	\Lambda_c(t)^2 + \gamma \rho_c(\gamma-p_1), &
				\Lambda_c(T)=0,
			\end{array}
			\\
			\begin{array}{ll}
				dY^c_t&=-\Big\{\frac{p_0+p_1 s_0-\gamma(s_0+\delta)}{2}+\frac{\rho_p(\gamma-p_1)}{2}q_t+\rho_p \gamma \rho_c \eta_c \lambda^2 (q_t - \mathbb E [q_t])+\frac{2}{k_c}\Big[K_c(t)(Y^c_t-\mathbb E [Y^c_t])\\
				&\quad + \Lambda_c(t)\mathbb E[Y^c_t]\Big]\Big\}dt + Z^{c,B}_tdB_t+Z^{c,W}_tdW_t,\\
				Y^c_T&=\frac{\lambda \gamma \rho_c}{2},
			\end{array}
			\\
			\begin{array}{ll}
				R_c'(t)  =& \eta_c \lambda^2 \rho_p^2 \mathbb V [q_t] -\frac{2}{k_c}(\mathbb V [Y^c_t]+\mathbb E [Y^c_t]^2)-\frac{2}{\ell_c-2K_c(t)}[\mathbb V[Z^{c,B}_t] +(\mathbb E [Z^{c,B}_t]+\frac{\ell_c \sigma_c}{2})^2], \\
				R_c(T)=&-\lambda \rho_p \mathbb E [q_T].
			\end{array}
		\end{array}\right.
	\end{split}
	\end{equation}
	So, we have 
	\[
		\widetilde{J}_p^\lambda(\mathbf{B}_p(\beta^c); \beta^c)=V_p^\lambda(\beta^c)
		\text{ and } 
		\widetilde{J}_c^\lambda(\mathbf{B}_c(\beta^p); \beta^p)=V_c^\lambda(\beta^p).
	\]
	Moreover, we have an explicit expression for the Nash equilibrium values which are given by
	\begin{equation}
		\begin{split}
			&V_p^\lambda(\beta^c)=  \Lambda_p(0) q_0^2 + 2 \mathbb E[Y^p_0] q_0 +R_p(0) 
			\quad \text{ and }\\&
			V_c^\lambda(\beta^p)= \Lambda_c(0) c_0^2 + 2 \mathbb E[Y^c_0] c_0 +R_c(0)  .
		\end{split}
	\end{equation}
\end{proposition}

\begin{remark}\label{K_and_Lambda}
	Notice that the first two equations in the systems \eqref{Sys_p} and \eqref{Sys_c} are one-dimensional Riccati differential equations, for which it is known that there exists a unique global solution given by Equation \eqref{Ks_Lambdas}.
	In the following we face more complicated Riccati equations (non-symmetric matrix Riccati equations) for which existence of solutions is not guaranteed.
	The fact that $K_p(t)$ and $K_c(t)$ are given by a hyperbolic tangent with a positive argument multiplied by a negative constant yields that $K_p(t)\leq 0$ and $K_c(t)\leq 0$, granting that $S_p(t)=N_p+e_2 K_p(t)e_2^\top$ and $S_c(t)=N_c+e_2 K_c(t)e_2^\top$ are negative definite for all $t \in [0,T],$ hence matching the assumptions made at the beginning of Sub-step 1.3. 
\end{remark}

\begin{proof}
	
To prove  the proposition we need to apply  Theorem \ref{Ver_thm}. So, let us check that its hypotheses are fulfilled.
	Fix a couple of strategies $\beta^p, \beta^c \in \mathcal{A}$.	
	First of all, condition i) is a consequence of the terminal conditions of systems \eqref{Sys_p} and \eqref{Sys_c}.
	Furthermore, we notice that assumption ii) is verified, for any $ \alpha^p \in \mathcal{A}$ (resp. for any $ \alpha^c \in \mathcal{A}$), because the fact that the processes $(K_p,\Lambda_p,Y^p,R_p) $ and $(K_c,\Lambda_c,Y^c,R_c)$ solve the systems  \eqref{Sys_p} and \eqref{Sys_c} yields that $\frac{d}{dt}\mathbb E [\mathcal{S}^{p,\alpha^p}_t]$ and $ \frac{d}{dt}\mathbb E [\mathcal S^{c,\alpha^c}_t]$ are negative and so the monotonicity of the functions $[0,T] \ni t \mapsto E [\mathcal{S}^{p,\alpha^p}_t] (\text{resp. }\mathbb E [\mathcal S^{c,\alpha^c}_t])$.
	Then,  by \eqref{opt_cont_cond}, we notice that, given $\beta^c \in \mathcal A$, $\frac{d}{dt}\mathbb E [\mathcal{S}^{p,\alpha^p}_t]=0, \text{ for all } t \in [0,T],$ if and only if, for all $t \in [0,T]$, we have 
	\begin{displaymath}
		\alpha^p_t
		=\eta^p_t
		= -S_p(t)^{-1}\Big[U_p(t)(q_t- \mathbb E[q_t])-V_p(t)\mathbb E[q_t]-(\xi^p_t-\bar{\xi}^p_t)- O_p(t)\Big], \quad \mathbb P\text{-a.s.},
	\end{displaymath}
	and analogously, given $\beta^p \in \mathcal A$, $0=\frac{d}{dt}\mathbb E [\mathcal S^{c,\alpha^c}_t], \text{ for all } t \in [0,T],$  if and only if, for all $t \in [0,T]$, we have 
	\begin{displaymath}
		\alpha^c_t
		=\eta^c_t
		=-S_c(t)^{-1}\left[U_c(t)(c_t- \mathbb E[c_t])-V_c(t)\mathbb E[c_t]-(\xi^c_t-\bar{\xi}^c_t)- O_c(t)\right], \quad \mathbb P\text{-a.s.}
	\end{displaymath}
	Hence, the strategies in \eqref{Best_r} satisfy iii) as well.\\
	Finally, let us check the admissibility of the strategies $\mathbf{B}_p(\beta^c)$ and $\mathbf{B}_c(\beta^p)$, i.e. $ \mathbf{B}_p(\beta^c)\in \mathcal{A}$ and $\mathbf{B}_c(\beta^p) \in \mathcal{A}$.
	We need to verify their square-integrability. Let us check it for $\mathbf{B}_p(\beta^c)$, the same can be done for $\mathbf{B}_c(\beta^p).$ The state variable $q=\{q_t\}_{t \in [0,T]}=\{q^{\mathbf{B}_p(\beta^c)}(t)\}_{t \in [0,T]}$ is the solution of a linear SDE and so it  satisfies $\mathbb E[\sup_{t \in [0,T]}|q_t|^2]<\infty$. 
	Furthermore, $S_p, U_p, V_p$, defined in \eqref{S_p_def}, are bounded, being continuous matrix-valued functions over a finite time-interval, and the process $(O_p,\xi^p)$ belongs to $L^2([0,T], \mathbb R^2)\times L^2_{\mathbb{F}}(\Omega \times [0,T], \mathbb R^2)$. This implies that the feedback control $\mathbf{B}_p(\beta^c) \in  L^2_{\mathbb F}(\Omega \times [0,T], \mathbb R^2) .$
\end{proof}

\subsection{Second step: well-posedness of the best response map system}

This subsection provides the proof of existence and uniqueness of solutions to the systems in \eqref{Sys_p} and \eqref{Sys_c},
\[
	K_p,K_c,\Lambda_p \text{ and }\Lambda_c \in L^{\infty}([0,T], \mathbb R_-), \quad
	R_p \text{ and } R_c \in L^{\infty}([0,T], \mathbb R),
\]
\[
	(Y^p,Z^{p,W}, Z^{p,B})\text{ and }(Y^c,Z^{c,W}, Z^{c,B}) \in S^2_{\mathbb F}(\Omega \times [0,T], \mathbb R)\times L^2_{\mathbb F}(\Omega \times [0,T], \mathbb R)^2,
\]
given the state controlled by the other player. \\

The fact that there exist unique $K_p,K_c,\Lambda_p,\Lambda_c \in L^{\infty}([0,T], \mathbb R_+)$ is straightforward (see Remark \ref{K_and_Lambda}). We also have explicit formulae  for them (see Equation \eqref{Ks_Lambdas}). Moreover the non-positivity of $ K_p$ and $K_c$ implies that the matrices $S_p$ and $S_c $, defined in \eqref{S_p_def}, are negative definite.\\
Now, consider the mean-field BSDE associated to the processes $(Y^p,Z^{p,W},Z^{p,B})$, given $K_p$ and $\Lambda_p:$
\[
	\left\{\begin{array}{l}
		dY^p_t =-\left\{\frac{s_0}{2}+\frac{\gamma \rho_c}{2} c_t +\rho_p \gamma \rho_c \eta_p \lambda^2 (c_t-\mathbb E [c_t])+\frac{2}{k_p}\left(K_p(t)(Y^p_t-\mathbb E [Y^p_t])+\Lambda_p(t)\mathbb E[Y^p_t]\right)\right\}dt\\
		\hspace{1.2cm} +Z^{p,B}_tdB_t+Z^{p,W}_tdW_t,\\
		Y^p_T=\frac{\lambda \rho_p}{2}.
	\end{array}\right.
\]
Existence and uniqueness of the solution $(Y^p,Z^{p,W},Z^{p,B})\in S^2_{\mathbb F}(\Omega \times [0,T], \mathbb R)\times L^2_{\mathbb F}(\Omega \times [0,T], \mathbb R)^2$ is a consequence of \cite[Theorem 2.1]{Li2017} and the fact that $c \in S^2_{\mathbb{F}}(\Omega \times [0,T], \mathbb{R})$ by the admissibility of the associated control $\beta^c$.\\
Finally, given $(K_p, \Lambda_p, (Y^p,Z^{p,W},Z^{p,B}))$, the linear ODE associated to $R_p$ in system \eqref{Sys_p} has a unique solution given by
\[
	\begin{split}
		R_p(t)=& - \lambda \gamma \rho_c \mathbb E [c_T]+\int_t^T \Bigg[-\eta_p \lambda^2 \gamma^2\rho_c^2 \mathbb V [c_u] +\frac{2}{k_p}\left(\mathbb V [Y^p_u]+\mathbb E [Y^p_u]^2\right)\\&
		+\frac{2}{\ell_p-2K_p(u)}\left(\mathbb V[Z^{p,W}_u]+\left(\mathbb E [Z^{p,W}_u]+\frac{\ell_p \sigma_p}{2}\right)^2\right)\Bigg]du.
	\end{split}
\]
The same arguments are used to prove existence and uniqueness for the processes 
$(Y^c,Z^{c,W}, Z^{c,B})$ in $S^2_{\mathbb F}(\Omega \times [0,T], \mathbb R)\times L^2_{\mathbb F}(\Omega \times [0,T], \mathbb R)^2$
and  
$ R_c \in L^{\infty}([0,T], \mathbb R) $. 
This ends the proof of existence and uniqueness for systems \eqref{Sys_p} and \eqref{Sys_c}.

\subsection{Third step: fixed point of the best response map}

Here, we prove the existence of a fixed point of the best response maps in order to get a Nash equilibrium.
First of all, for convenience of notation, we rewrite the two-dimensional state variable as $X_t:=(	q_t, c_t)^\top,$ for all  $t \in [0,T]$, and so its linear dynamics is given by the following SDE
\begin{align}
	&dX_t 
	=\begin{pmatrix}
		dq_t  \\
		dc_t
	\end{pmatrix} 
	=
	\begin{pmatrix}
		u_t  \\
		v_t
	\end{pmatrix}dt
	+
	\begin{pmatrix}
		z_t  \\
		0
	\end{pmatrix} dW_t
	+
	\begin{pmatrix}
		0  \\
		y_t
	\end{pmatrix} dB_t,
\end{align}
with a deterministic initial condition $X_0 = (q_0, c_0)^\top \in \mathbb{R}_+^2.$
Then, we have
\begin{align*}
	dX_t=b \alpha_t dt + \sigma^W \alpha_t dW_t + \sigma^B \alpha_t dB_t,
\end{align*}
with
\begin{align*}
\begin{split}
	&b =
	\begin{pmatrix}
		1 & 0 & 0 & 0 \\
		0 & 0 & 1 & 0
	\end{pmatrix},\quad
	\sigma^W =
	\begin{pmatrix}
		0 & 1 & 0 & 0\\
		0 & 0 & 0 & 0
	\end{pmatrix},
	\quad
	\sigma^B =
	\begin{pmatrix}
		0 & 0 & 0 & 0 \\
		0 & 0 & 0 & 1
	\end{pmatrix}.
	\end{split}
\end{align*}
Then, we rewrite explicitly the form that a candidate equilibrium feedback control $\alpha^*=((\alpha^{*,P})^\top,(\alpha^{*,C})^\top)^\top$ should have, together with the backward dynamics of the corresponding process $Y=((Y^p)^\top,(Y^c)^\top)^\top$  (we write $Z^W$ for $((Z^{p,W})^\top,(Z^{c,W})^\top)^\top$, respectively $Z^B$ for $((Z^{p,B})^\top,(Z^{c,B})^\top)^\top$),\footnote{Here, we have omitted, in all the processes but $\alpha$, the superscript $*$ in order to have a simpler notation.}

\begin{align}\label{din_alpha}
	&\alpha^*_t- \bar \alpha^*_t=\Delta(t) \left(X_t-\bar X_t\right)+\Gamma \left(Y_t-\bar Y_t\right)+H^W(t)\left(Z_t^W-\bar Z_t^W\right)+H^B(t)\left(Z_t^B-\bar Z_t^B\right),\\&\nonumber
	\bar \alpha^*_t=\widehat \Delta(t)\bar X_t+ \Gamma\bar Y_t+ H^W(t)\bar Z_t^W+ H^B(t) \bar Z_t^B+ \Theta(t),
\end{align}
\begin{align}\label{din_Y}
	dY_t=\left[\Xi \left(X_t-\bar X_t\right)+ \Phi(t) \left(Y_t-\bar Y_t \right)\right]dt + \left[ \widehat \Xi \bar X_t +\widehat \Phi(t)\bar Y_t+ \Psi \right]dt+Z_t^BdB_t+Z^W_t dW_t,
\end{align}
with
\begin{align*}
	&\Delta(t)=
	\begin{pmatrix}
		\frac{2}{k_p}K_p(t) & 0 \\
		0 & 0\\
		0 &\frac{2}{k_c}K_c(t)  \\
		0 & 0
	\end{pmatrix}, \quad
	\widehat \Delta(t) =
	\begin{pmatrix}
		\frac{2}{k_p}\Lambda_p(t) & 0 \\
		0 & 0\\
		0 &\frac{2}{k_c}\Lambda_c(t)  \\
		0 & 0
	\end{pmatrix},\quad
	\Gamma
	=
	\begin{pmatrix}
		\frac{2}{k_p} & 0 \\
		0 & 0\\
		0 &\frac{2}{k_c} \\
		0 & 0
	\end{pmatrix}, 
\end{align*}
\begin{align}
	\nonumber
	&\Theta(t)=
	\begin{pmatrix}
		0 \\
		\sigma_p(1-2\frac{K_p(t)}{\ell_p})^{-1}\\
		0 \\
		\sigma_c(1-2\frac{K_c(t)}{\ell_c})^{-1}
	\end{pmatrix}, \quad
	H^W(t)=
	\begin{pmatrix}
		0  & 0\\
		 \frac{2}{\ell_p -2 K_p(t)} &0\\
		0  & 0\\
		0  & 0 
	\end{pmatrix}, \quad
	H^B(t)=
	\begin{pmatrix}
		0  & 0\\
		0  & 0\\
		0  & 0\\
		0 &  \frac{2}{\ell_c -2 K_c(t)}
	\end{pmatrix},
\end{align}
 and 
 $\Xi,$ $\widehat{\Xi},$ $\Phi(t)$, $\widehat{\Phi}(t)$ and $\Psi$ as defined at the beginning of Section \ref{sect_main_res}.
\\

Now, as an ansatz for $Y$, we assume $Y$ linear in the state:
\begin{equation}\label{ansatz_Y}
	Y_t=\pi(t)(X_t-\bar X_t)+\widehat \pi(t) \bar X_t+\zeta_t, 
\end{equation}
with $\pi,\widehat \pi$ deterministic $\mathbb R^{2\times 2}$-valued processes and $\zeta \in S^2_{\mathbb F}(\Omega \times [0,T], \mathbb R^2)$ satisfying the SDE
\begin{align}\label{zeta}
	d\zeta_t=\psi_t dt+\phi_t^W dW_t+\phi_t^B dB_t,
	\quad 
	\zeta_T=\frac{1}{2}\lambda
	(\rho_p,\gamma \rho_c)^\top,
\end{align}
for some $\psi, \phi^B, \phi^W$ in suitable spaces.
The affine term in the expression \eqref{ansatz_Y} allows $Y$ to have some extra stochasticity apart from the linear dependency on the state. Furthermore, the terminal condition in \eqref{zeta} guarantees that $Y$ satisfies its terminal condition.\\
An application of It\^o's formula to the ansatz \eqref{ansatz_Y} yields 
\begin{equation}\label{din_Y_ans}
	\begin{split}
		dY_t=&[\pi'(t)(X_t-\bar X_t)+\pi(t)b(\alpha_t^*-\bar\alpha_t^*)+\psi_t-\bar\psi_t]dt
		+(\widehat\pi'(t)\bar X_t+\widehat \pi(t)b\bar\alpha_t^*+ \bar\psi_t)dt\\&
		+(\pi(t)\sigma^W\alpha^*_t+\phi^W_t)dW_t
		+(\pi(t)\sigma^B\alpha^*_t+\phi^B_t)dB_t.
	\end{split}
\end{equation}
If we match the two dynamics of $Y$ in Equations \eqref{din_Y} and \eqref{din_Y_ans}, and then replace $Y$ with its ansatz \eqref{ansatz_Y} and $\alpha^*$ with its feedback form \eqref{din_alpha}, we get  the following system of equations:

\begin{equation}\label{Eq_big_sys_Z}
	\begin{split}
	&\left\{\begin{array}{ll}
		\pi'(t)(X_t-\bar X_t)+\pi(t)b(\bold I -H^W(t) \pi(t)\sigma^W-	H^B(t) \pi(t)\sigma^B)^{-1}[(\Delta(t)\\
		\quad +\Gamma\pi(t))(X_t-\bar X_t)+\Gamma(\zeta_t-\bar\zeta_t)+H^W(t)(\phi_t^W-\bar\phi_t^W)+H^B(t)(\phi_t^B-\bar \phi_t^B)]\\
		\quad +\psi_t-\bar\psi_t = \Xi(X_t-\bar X_t)+ \Phi(t) (\pi(t)(X_t-\bar X_t)+\zeta_t-\bar 	\zeta_t)
		\\
		\\
		\widehat\pi'(t)\bar X_t+\widehat \pi(t)b(\bold I -H^W(t) \pi(t)	\sigma^W-H^B(t) \pi(t)\sigma^B)^{-1}[(\widehat\Delta(t)+\Gamma\widehat	\pi(t))\bar X_t
		+\Gamma\bar\zeta_t\\
		\quad +\Theta(t) + H^W(t)\bar\phi_t^W+H^B(t)\bar\phi_t^B]
		+ \bar\psi_t =
		\widehat \Xi \bar X_t +\widehat \Phi(t)(\widehat \pi(t) \bar X_t + 	\bar \zeta_t)+ \Psi
		\\ 
		\\
		Z_t^B=\pi(t)\sigma^W\alpha^*_t+\phi^W_t
		\\
		Z^W_t=\pi(t)\sigma^B\alpha^*_t+\phi^B_t.
	\end{array}\right.
	\end{split}
\end{equation}
Finally, exploiting the fact that:
\begin{equation}
	b\left(\bold I -H^W(t) \pi(t)\sigma^W- H^B(t) \pi(t)\sigma^B\right)^{-1}=b,
\end{equation}
we find the equations that the coefficients $(\pi,\widehat \pi, \psi, \phi^W,\phi^B)$ in the ansatz for $Y$ should solve in order to provide a fixed point of the best response map:

\begin{equation}\label{Sys_Nash}
	\begin{split}
		&\left\{\begin{array}{l}
			\begin{array}{ll}
				\pi'(t)=\Xi+\Phi(t) \pi(t)+\pi(t)\Phi(t)+\pi(t)R\pi(t),&
				\pi(T)=0,
			\end{array}
			\\
			\begin{array}{ll}
				\widehat\pi'(t)=\widehat\Xi+\widehat\Phi(t) \widehat\pi(t)+	\widehat\pi(t)\widehat\Phi(t)+\widehat\pi(t)R\widehat\pi(t),&
				\widehat\pi(T)=0,
			\end{array}
			\\
			\begin{array}{ll}
				d\zeta_t=\psi_t dt+\phi_t^W dW_t+\phi_t^B dB_t,&
				\zeta_T=\frac{1}{2}\lambda
				(\rho_p, \gamma \rho_c)^\top,
			\end{array}
			\\
			\psi_t=\psi_t-\bar\psi_t+\bar\psi_t=\Big(\pi(t)R+\Phi(t)\Big)(\zeta_t-\bar\zeta_t)+\left(\widehat\pi(t)R+\widehat\Phi(t)\right)\bar\zeta_t+\Psi,
		\end{array}\right.
	\end{split}
\end{equation}
where $
	R=
	\begin{pmatrix}
		- 2/k_p & 0 \\
		0 & -2/k_c 
	\end{pmatrix}.
$ 
In fact, inserting $(Y,Z)$ from the ansatz and Equation \eqref{Eq_big_sys_Z} into the best response given by Equations \eqref{din_alpha} provides an equilibrium strategy $\alpha^*$ in feedback form which is computed in details in the next step. 

\begin{remark}
To obtain explicit expressions for $\alpha^*$ and  $Z$, we have used Assumption \textbf{(A2)} in Theorem \ref{thm_Nash}. Indeed, such a condition is needed for the invertibility of the matrices $D(t):=(\bold I -H^W(t) \pi(t)\sigma^W-H^B(t) \pi(t)\sigma^B)$, $t \in [0,T]$, that appear in 
\begin{align} \label{eq:ZwZb}
	&Z^W_t
	=\phi_t^W+\pi(t)\sigma^W\alpha^*_t,
	\qquad 
	Z^B_t
	=\phi_t^B+\pi(t)\sigma^B\alpha^*_t,
\end{align}
where 
\begin{equation*}
	\alpha^*_t= D(t)^{-1}[(\Delta(t)+\Gamma\pi(t))(X_t-\bar X_t) + 	(\widehat\Delta(t)+\Gamma\widehat \pi(t))\bar X_t+\Gamma \zeta_t 	+H^W(t)\phi_t^W+H^B(t)\phi_t^B + \Theta(t)].
\end{equation*}
\end{remark}

\subsection{Fourth step: Nash equilibrium strategies}
In order to complete the proof of the main theorem, we are left with showing that  the system \eqref{Sys_Nash} has a unique solution over the finite time interval $[0,T]$.   
The equations associated to $t \mapsto (\pi(t), \widehat \pi(t))$ are non-symmetric matrix Riccati equations for which there is no general condition ensuring the global existence of solutions. Nevertheless, the regularity of the coefficients and the Picard-Lindel\"{o}f Theorem ensure the local existence and uniqueness of solutions over a compact interval $[0,T_{max}]$\footnote{Despite it is not possible in general to obtain an explicit characterization of $T_{max}$, we notice that we did not observe any explosion for all typical values of the parameters we have considered in the numerical experiments (ref. Section \ref{sec:numeric}). 
}. Thus, we recover the existence and uniqueness condition in Assumption \textbf{(A1)} of Theorem \ref{thm_Nash} choosing a time horizon $T$ small enough, namely $T < T_{max}$.
Then, for a given $(\pi,\widehat \pi)$, the process $(\zeta, \psi, \phi^W,\phi^B)$ evolves according to the following linear mean field BSDE:
\begin{equation}\label{sys_zeta}
	\begin{split}
		&d\zeta_t=\psi_t dt+\phi^W_t dW_t+\phi^B_t dB_t, 
		\quad
		\zeta_T=\frac{1}{2}\lambda (\rho_p, \gamma \rho_c)^\top,
		\\&
		\psi_t=\psi_t-\bar\psi_t+\bar\psi_t=(\pi(t)R+\Phi(t))(\zeta_t-\bar\zeta_t)+(\widehat\pi(t)R+\widehat\Phi(t))\bar\zeta_t+\Psi.
	\end{split}
\end{equation}
Exploiting once more \cite[Theorem 2.1]{Li2017}, we have a unique solution $(\zeta,\phi^W, \phi^B)  \in S^2_{\mathbb F}(\Omega \times [0,T], \mathbb R^2)\times L^2_{\mathbb F}(\Omega \times [0,T], \mathbb R^2)^2.$
Furthermore, we notice that the drift $\psi$ in the system \eqref{sys_zeta} does not depend on $\phi^W$ and $\phi^B$ and all the coefficients involved in the second line of  \eqref{sys_zeta} are deterministic. Moreover, the terminal condition is also deterministic. Thus, the unique solution $(\zeta, \phi^W, \phi^B)$ to this system is given by $(h, 0, 0)$, where $h:[0,T] \to \mathbb{R}^2$ is the unique (deterministic) solution to the following backward linear ODE: 
\begin{equation}\label{sys_h}
	\begin{split}
		&\left\{\begin{array}{l}
			dh(t)=\left\{\left[\widehat\pi(t)R+\widehat\Phi(t)\right]h(t)+\Psi\right\}dt,\\
			h(T)=\frac{1}{2}\lambda
			(\rho_p,\gamma \rho_c)^\top.
		\end{array}\right.
	\end{split}
\end{equation}
So, the system of ODEs and SDEs in \eqref{Sys_Nash} reduces to the one made up of Equations \eqref{Sys_Nash_simple_1} and \eqref{Sys_Nash_simple_2}.

We write the  Nash equilibrium strategies $\alpha^*=((\alpha^{*,P}),(\alpha^{*,C}))^\top = ((u^*,z^*)^{\top},(v^*,y^*)^\top)^\top$ explicitly as
	\begin{equation}
		\alpha^*_t
		= D(t)^{-1}(\Delta(t)+\Gamma \pi(t))(X_t-\bar{X}_t)+D(t)^{-1}(\widehat\Delta(t)+\Gamma \widehat\pi(t))\bar{X}_t + D(t)^{-1}(\Gamma h(t)+ \Theta(t)),
	\end{equation}
	that is
	\begin{align}
		\nonumber
		&u^*_t = \frac{2}{k_p}\left[(K_p(t)+\pi_{11}(t))(q_t-\bar q_t)+\pi_{12}(t)(c_t-\bar c_t)+(\Lambda_p(t)+\widehat \pi_{11}(t))\bar q_t+\widehat\pi_{12}(t)\bar c_t +h_1(t)\right],\\&\nonumber
		z^*(t) = \frac{\sigma_p \ell_p}{\ell_p-2(K_p(t)+\pi_{11}(t))},\\&\nonumber
		v^*_t = \frac{2}{k_c}\left[(K_c(t)+\pi_{22}(t))(c_t-\bar c_t)+\pi_{21}(t)(q_t-\bar q_t)+(\Lambda_c(t)+\widehat \pi_{22}(t))\bar c_t+\widehat\pi_{21}(t)\bar q_t +h_2(t)\right],\\&\nonumber
		y^*(t) = \frac{\sigma_c \ell_c}{\ell_c-2(K_c(t)+\pi_{22}(t))},
	\end{align}
where $K_p, K_c, \Lambda_p, \Lambda_c$ are defined in \eqref{Ks_Lambdas} and $\pi, \widehat{\pi}$ and $h$ are respectively the solutions to the systems \eqref{Sys_Nash_simple_1}, \eqref{Sys_Nash_simple_2}.

Finally, we derive the corresponding equilibrium dynamics for the state

\[
\begin{split}
	dX_t  
	= &
	\Bigg\{
	\begin{pmatrix}
		\frac{2}{k_p}(K_p(t)+\pi_{11}(t)) & \frac{2}{k_p}\pi_{12}(t) \\
		\frac{2}{k_c}\pi_{21}(t) & \frac{2}{k_c}(K_c(t)+\pi_{22}(t))
	\end{pmatrix} (X_t - \bar{X}_t)\\&
	+ 
	\begin{pmatrix}
		\frac{2}{k_p}(\Lambda_p(t)+\widehat\pi_{11}(t)) & \frac{2}{k_p}\widehat\pi_{12}(t) \\
		\frac{2}{k_c}\widehat\pi_{21}(t) & \frac{2}{k_c}(\Lambda_c(t)+\widehat\pi_{22}(t))
	\end{pmatrix} \bar{X}_t
	+
	\begin{pmatrix}
		\frac{2}{k_p}h_1(t) \\
		\frac{2}{k_c}h_2(t)
	\end{pmatrix}
	\Bigg\} dt\\&
	+ 
	\begin{pmatrix}
		\frac{\sigma_p \ell_p}{\ell_p-2(K_p(t)+\pi_{11}(t))} \\
		0
	\end{pmatrix}
	dW_t
	+
	\begin{pmatrix}
		0 \\
		\frac{\sigma_c \ell_c}{\ell_c-2(K_c(t)+\pi_{22}(t))}
	\end{pmatrix}
	dB_t,  \quad t \in [0,T],
\end{split}
\]
which is a linear mean-field SDE, hence admitting a unique solution.

\section{Numerics}\label{sec:numeric}
\label{sec:riskpremia}

We consider the following parameters setting $T=1$, $k_p=$ $k_c=5$,  $\sigma_p=$ $\sigma_c=10$,  $q_0=$ $c_0=100$,   $s_0=50$, $\rho_p=$ $\gamma \rho_c=0.5$ and $\gamma=1.2$, $ \delta=5$, $p_0=$ $2 s_0+\gamma \delta$, and $p_1=\gamma-1$. With this parametrisation, the players are symmetric in the sense that they have the same absolute effect on the price and they share the same costs of average production rate or consumption rate.  Moreover, if they shared the same risk aversion parameters ($\eta_p=\eta_c$) and the same costs of volatility control ($\ell_p = \ell_c$), then the strategies of the producer $(u^\ast,z^\ast)$ and of the consumer $(v^\ast,y^\ast)$ would be identical. The initial condition have been chosen to be close to a long-run stationary equilibrium that we observe when we take large $T$ and thus, allowing to avoid potential transitory effects. 

In the next sub-sections, we illustrate first the effect of the risk aversion parameters on the forward agreement indifference price when every other parameter is fixed. Second, we show how different combinations of risk aversions and volatility control costs can lead to the same forward agreement indifference price and volume.

\subsection{The effect of risk aversion}
\begin{figure}[thb!]
\centering
\begin{tabular}{c  c} 
(a) & (b) \\
\includegraphics[width=0.45\textwidth]{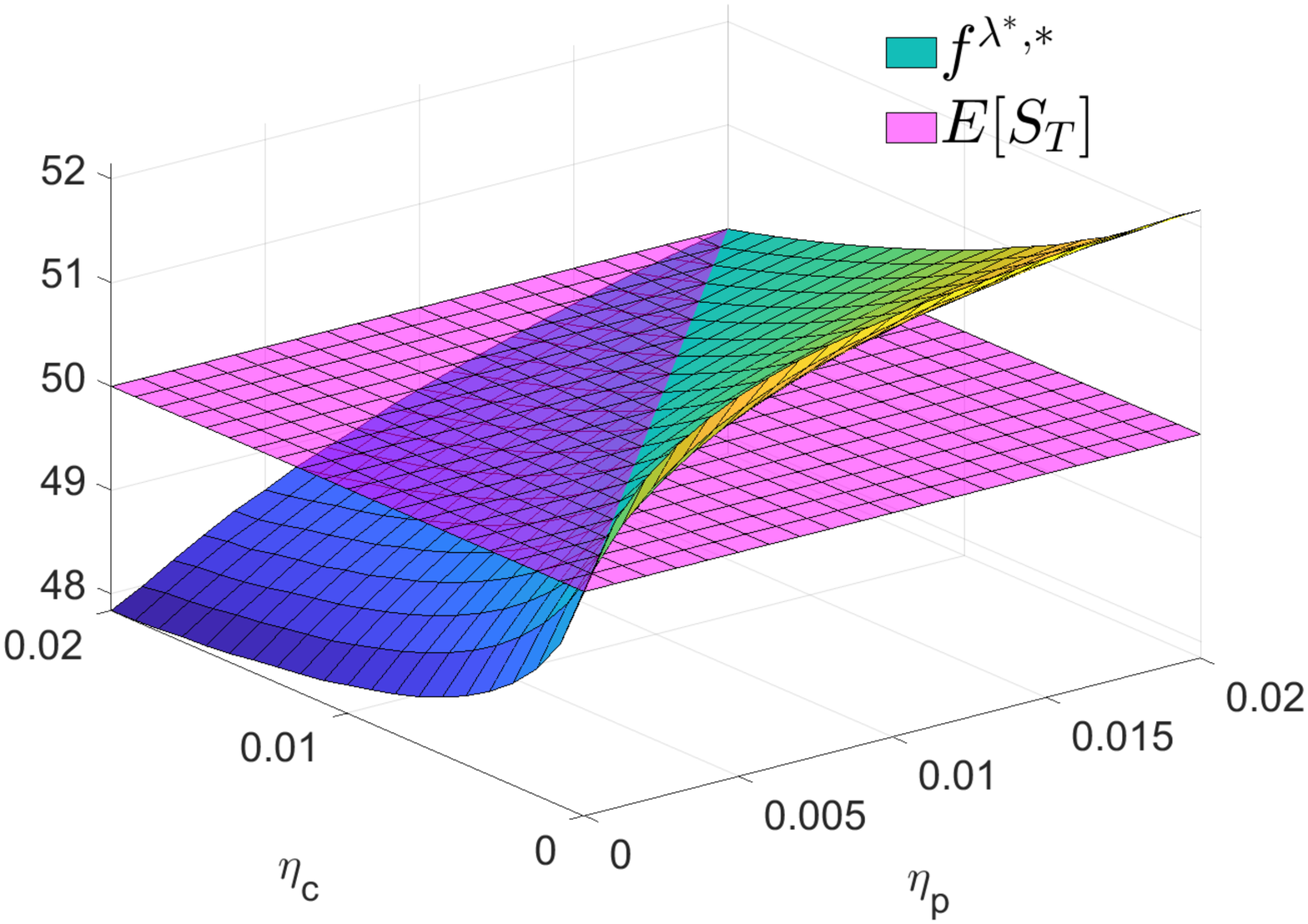} & \includegraphics[width=0.45\textwidth]{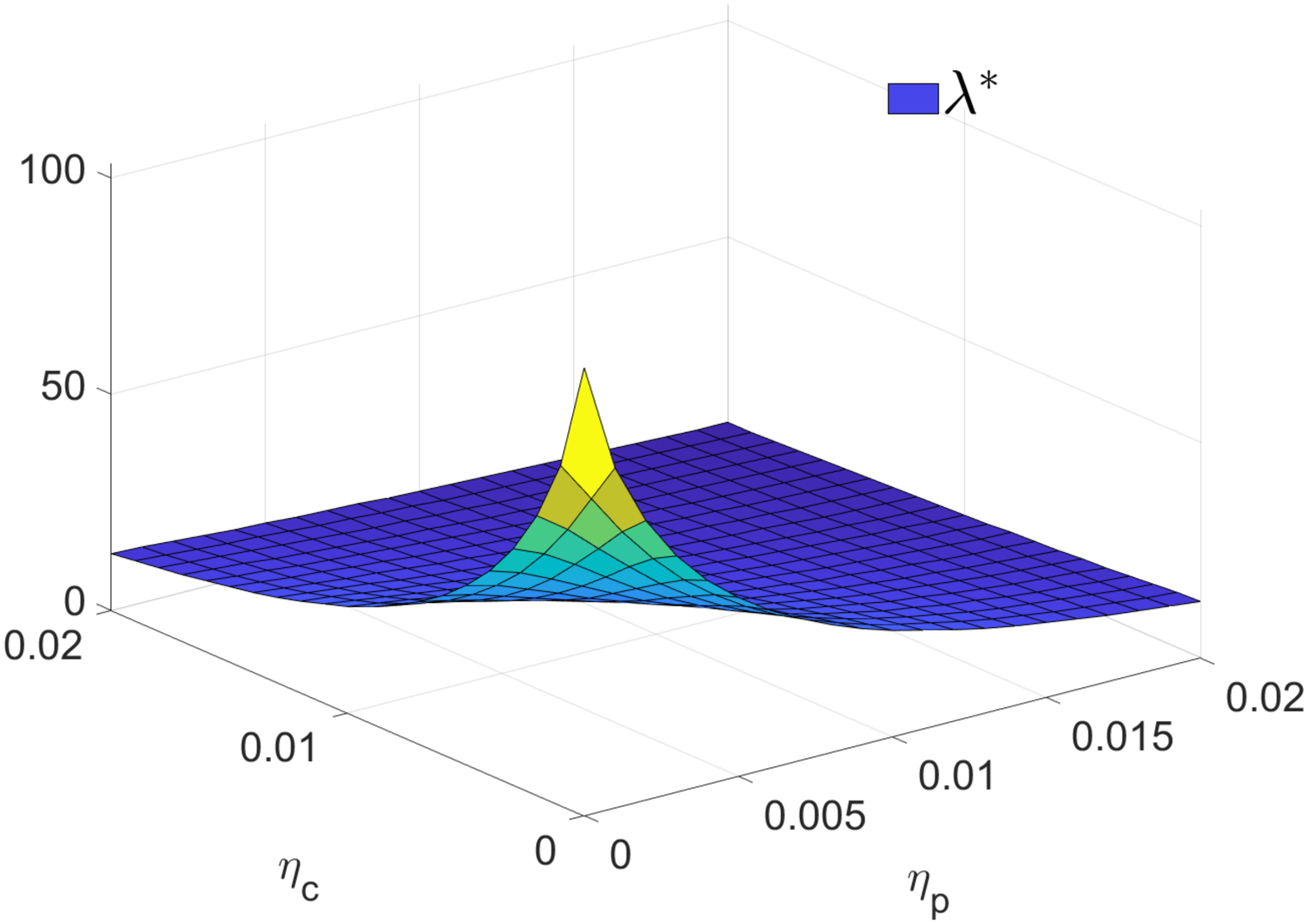} \\
(c) & (d) \\
\includegraphics[width=0.45\textwidth]{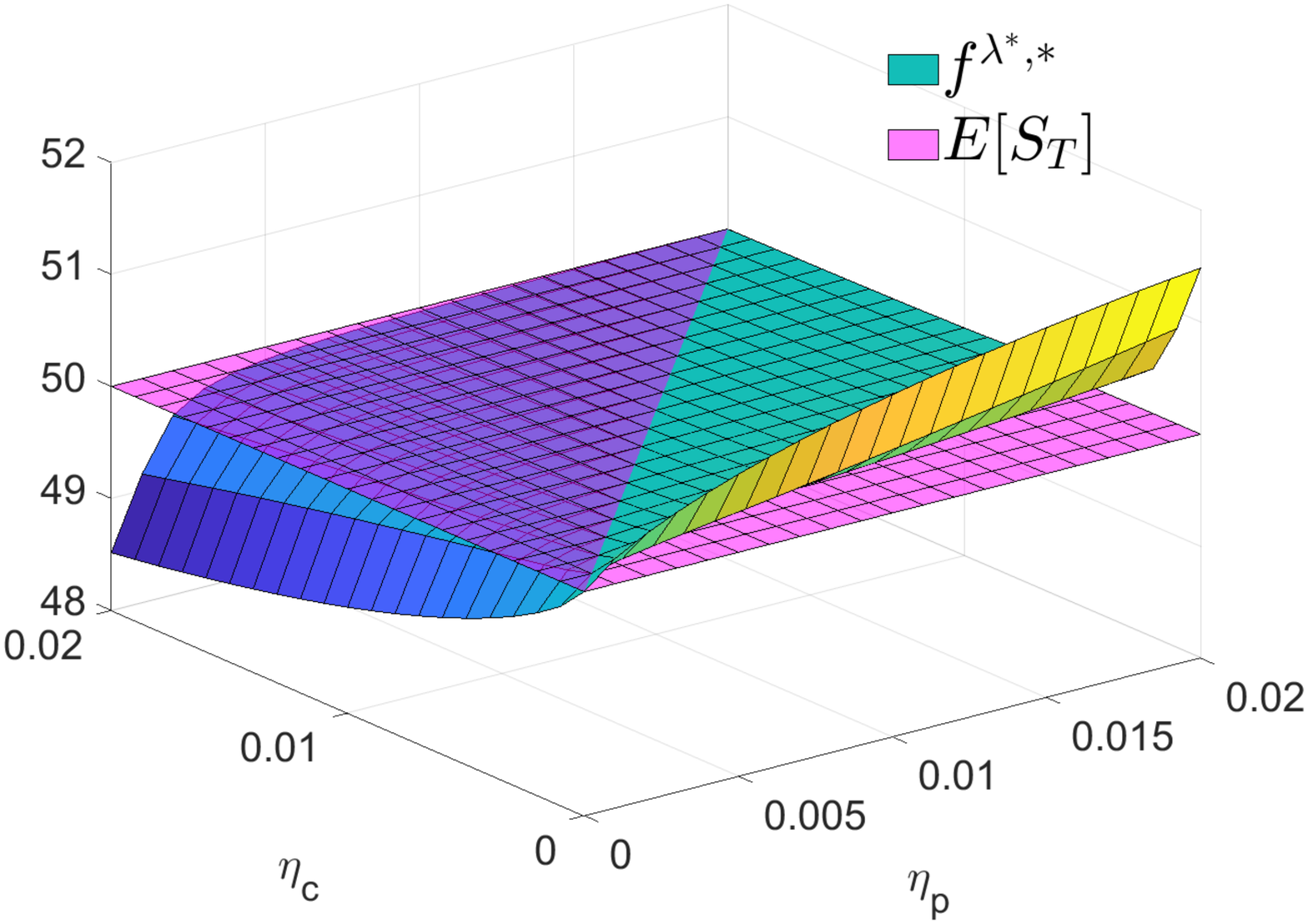} & \includegraphics[width=0.45\textwidth]{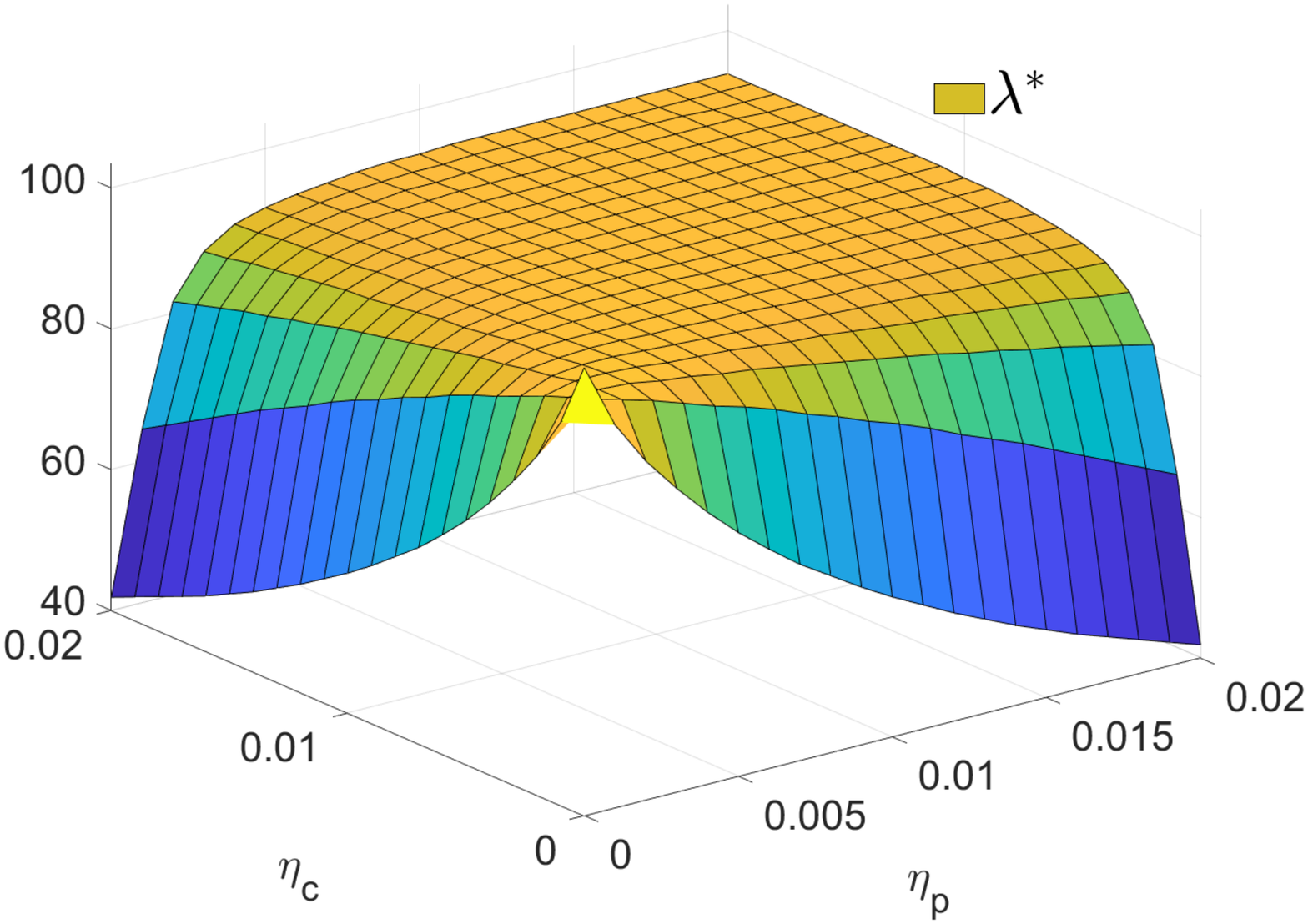}
\end{tabular}
\caption{ (a) and (b) $\ell_p=\ell_c=5$, (c) and (d) $\ell_p=\ell_c=0.7$. }\label{fig:riskprime}
\end{figure}

Figure~\ref{fig:riskprime}  presents the unitary forward agreement indifference price ${f}^{\lambda^\ast,\ast} := {F}^\ast_{{\lambda}^\ast}/\lambda^\ast$ and the volume that the players agreed upon when the costs of volatility control are high (Figure~\ref{fig:riskprime}~(a) and (b)) and when they are low (Figure~\ref{fig:riskprime}~(c) and (d)). We find that ${f}^{\lambda^\ast,\ast}$ is higher (resp. lower) than the expected spot price when the producer is more (resp. less) risk-averse than the consumer, which is consistent with both the economic intuition and the hedging pressure theory, once recalled that in our model players act as speculators on the forward market. In hedging pressure theory (see \cite{DNV00} and \cite{ELV19}), the risk premium is determined by the relations between risk aversions of producers, consumers, storers and speculators. It extends Keynes's normal backwardation theory which claims that in commodity markets, the forward price should be lower than the expected spot price because the producer would be ready to pay a premium to avoid being exposed to price risk on his production. In our case, the most risk-averse speculator obtains the appropriate premium to enter into the agreement. This property holds whatever the level of volatility control costs. We see on Figure~\ref{fig:riskprime} that the producer is requiring a positive premium to accept the  risk coming from his financial position. Regarding the exchanged volume, we observe that it can be both nonincreasing or nondecreasing in the risk aversion parameters of the players, depending on the costs of volatility control.
When the volatility manipulation costs are high for both players, there is a low trading volume even when both players have a high risk aversion. On the other side, when the volatility manipulation costs are low, there is a low trading volume when only one of the player has a high risk aversion but the trading volume is huge when both players have a high risk aversion. This could be explained by the fact that in the latter case the players can act on their volatilities (almost costlessly) to stabilize the spot price and hence they would be willing to trade more.

\subsection{Joint effect of risk aversion and volatility control cost}

\begin{figure}[thb!]
\centering
\begin{tabular}{c  c} 
(a) & (b) \\
\includegraphics[width=0.45\textwidth]{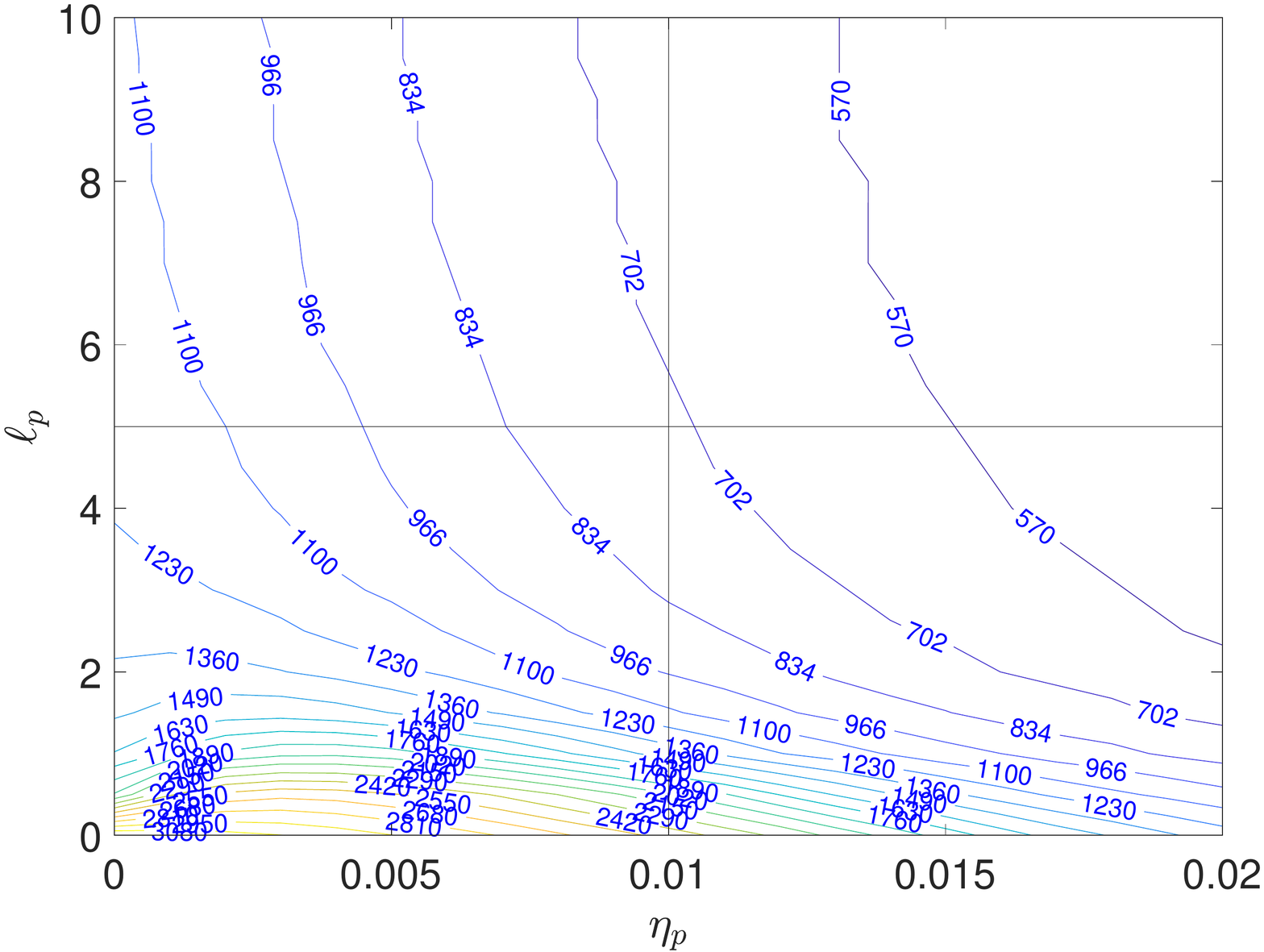} & \includegraphics[width=0.45\textwidth]{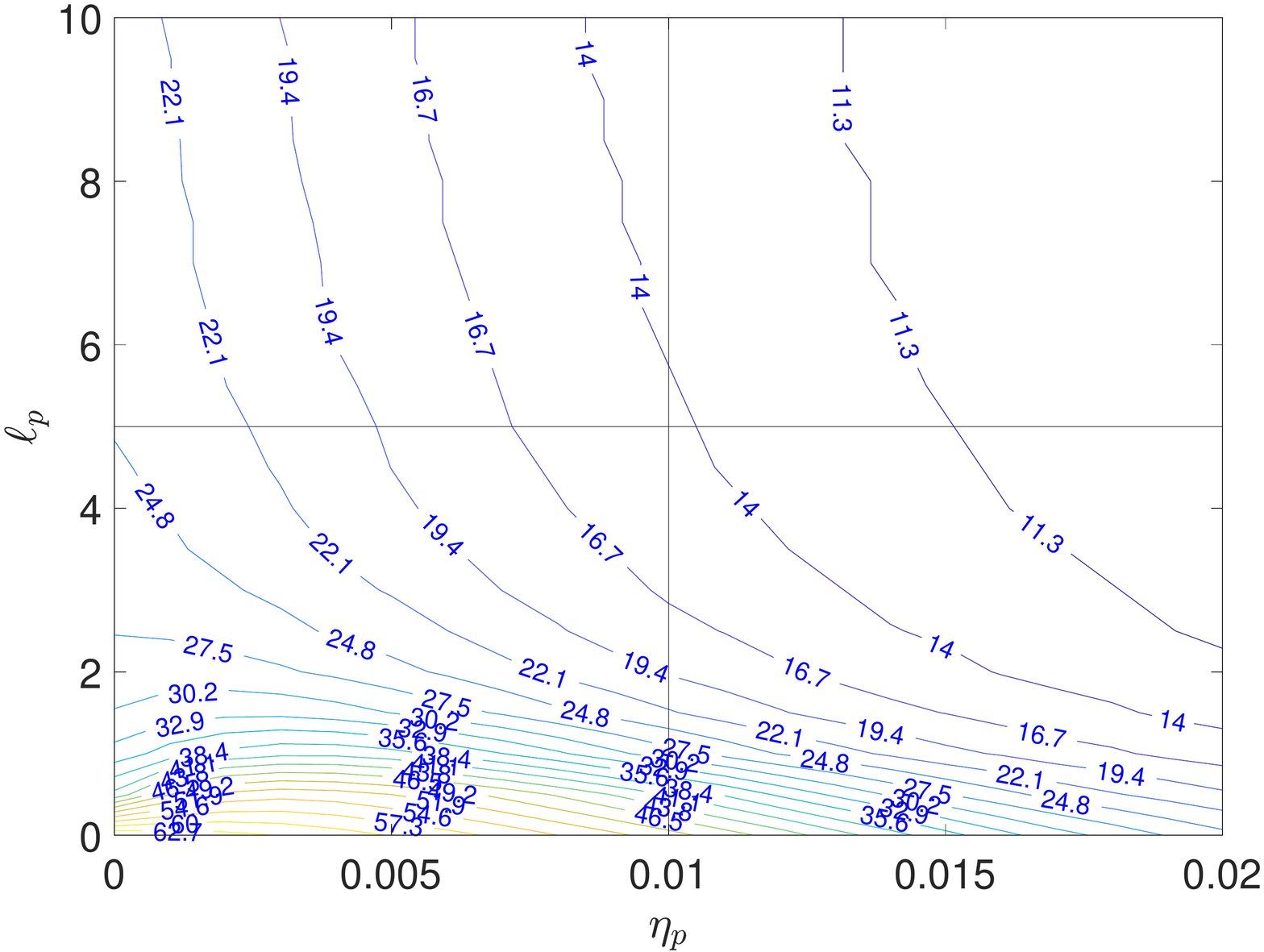} \\
(c) & (d) \\
\includegraphics[width=0.45\textwidth]{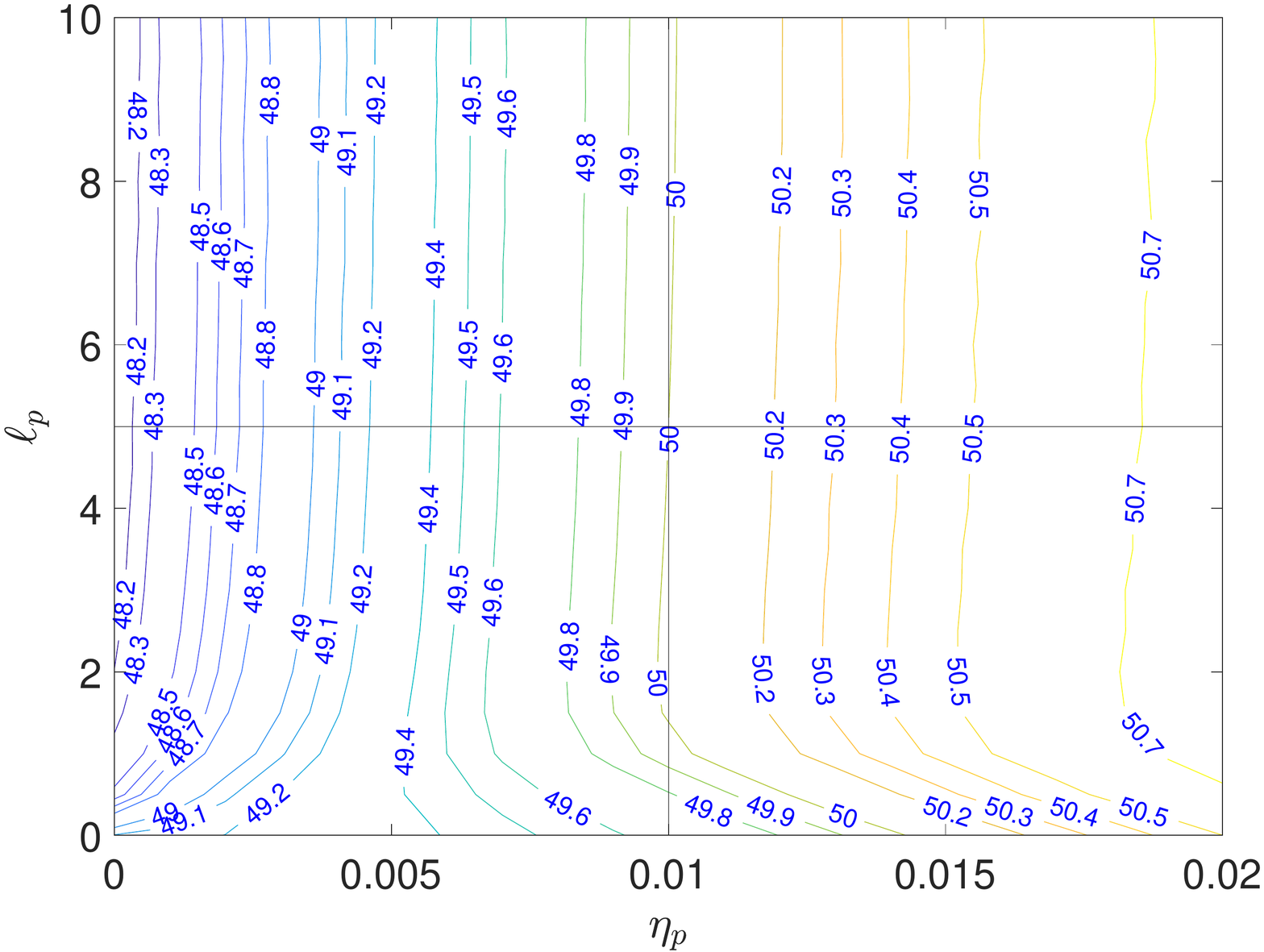} & \includegraphics[width=0.45\textwidth]{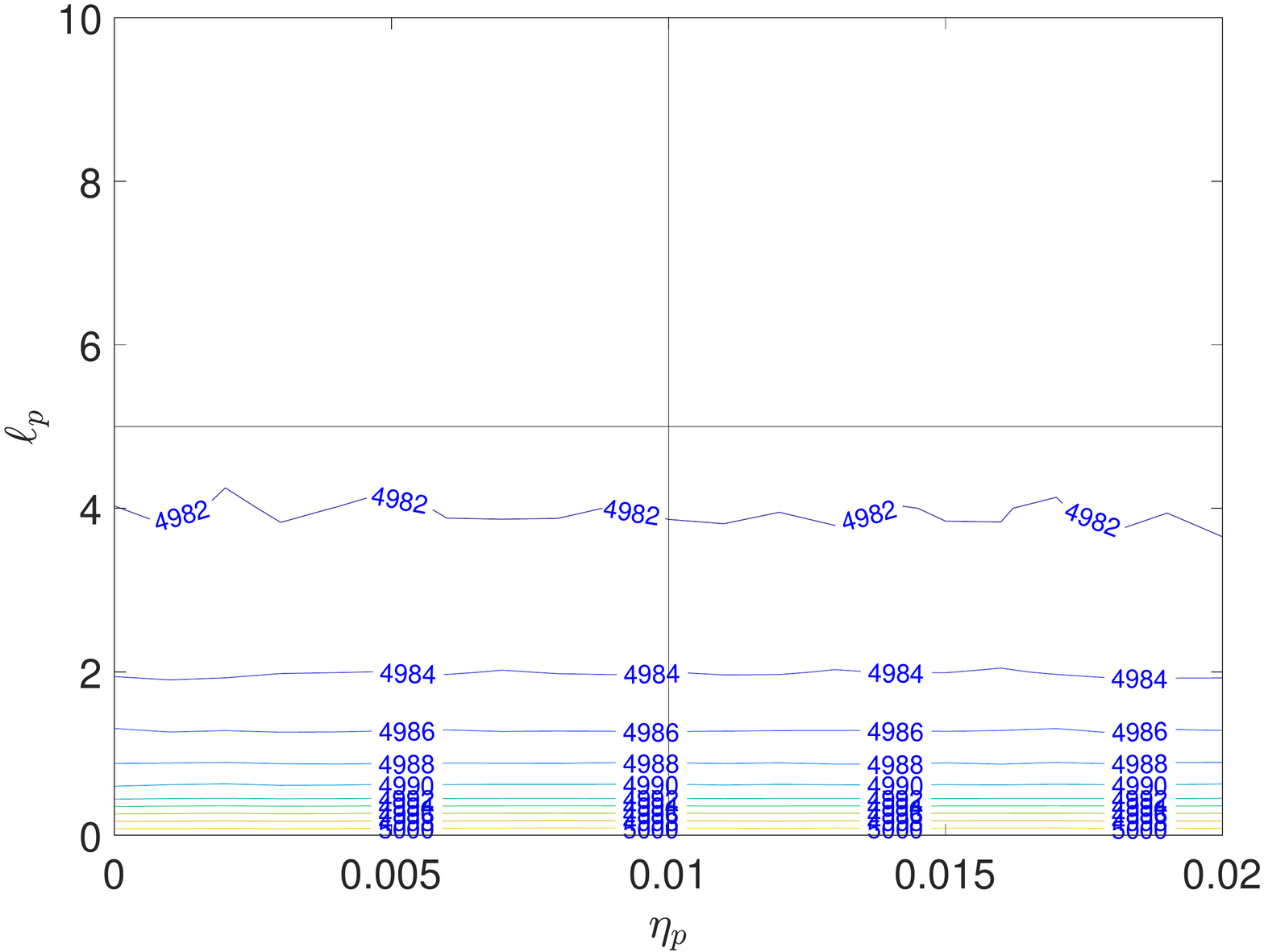}
\end{tabular}
\caption{ Level lines of (a) the forward agreement price ${F}^{\ast}_{\lambda^\ast}$, (b) the traded quantity $\lambda^\ast$, (c) the per unit agreement price ${f}^{\lambda^\ast,\ast}={F}^\ast_{{\lambda}^\ast}/\lambda^\ast$, (d) the value of the producer's equilibrium payoff $ J_p^\ast(\lambda^\ast,{F}^\ast_{{\lambda}^\ast})$. }\label{fig:substitute}
\end{figure}

We freeze now the risk aversion parameter and the cost for controlling the volatility of the consumer at $\eta_c=0.01$ and $\ell_c=5$, and observe the agreement price, the traded volume, the per unit agreement indifference price and the equilibrium payoff at the agreement of the producer. Results are provided in Figure~\ref{fig:substitute}, when the producer's risk aversion parameter $\eta_p$ and his volatility manipulation cost $\ell_p$  vary. The vertical and horizontal lines in each graph are set to the values of $\ell_c$ and $\eta_c$.

We observe a sort of  "substitution effect" between $\eta_p$ and $\ell_p$ in the sense that for a producer with a given combination of risk aversion and volatility control cost, we can find another producer \emph{trading at the same agreement price} with a higher risk aversion and a low volatility control cost (Figure~\ref{fig:substitute}~(a)). We observe that this phenomenon occurs also for the traded quantity (Figure~\ref{fig:substitute}~(b)). This substitution makes sense in our model where volatility represents a cost for the producer that can be mitigated either by requiring a payment to bear this volatility or by paying the cost to reduce it. 
We note that for a fixed value of $\eta_p$, the lower the value of $\ell_p$, the larger the forward agreement price {\em and} the traded volume. The Figure~\ref{fig:substitute}~(c) gives the resulting unitary agreement forward price. The volatility control cost has little effect on the per unit forward price compared to the risk aversion parameter. This figure is a way of showing that when the volatility control costs are high, the producer has little alternative than asking for a premium to enter in forward agreement, and thus, the price is basically determined by his risk-aversion. 

To conclude, we note that the producer's equilibrium payoff is independent of the value of $\eta_p$ (Figure~\ref{fig:substitute}~(d)) because, by definition of the agreement forward price, it is always equal to $J_p^\ast (0,0)$, which is independent of $\eta_p$.

\vspace{1cm}

\textbf{Funding:} this work was supported by Universit\`a degli Studi di Padova [Grant BIRD 190200/19]. 

\appendix

\section{Proof of Theorem \ref{Ver_thm} (Verification Theorem)}

For any $ \alpha^p \in \mathcal{A}$ (resp. for any $ \alpha^c \in \mathcal{A}$), the map $[0,T] \ni t \mapsto \mathbb{E}[\mathcal{S}^{p,{\alpha^p}}_t] \Big(\text{resp. } \mathbb{E}[\mathcal{S}^{c,{\alpha^p}}_t]\Big)$ is well-defined (it does not explode in finite time), because of the condition \eqref{growth_cond} and the linear structure of the SDEs for the state variables \eqref{eq_state_game}.\\
	Assumptions i) and ii) yields that: for any $ \alpha^p \in \mathcal{A}^p$,
	\[
		\begin{split}
			\mathbb{E}[&w_0^p(q_0,\bar q_0)]
			= \mathbb{E}[\mathcal{S}^{p,\alpha^p}_0]
			\overset{\text{ii)}}{\geq} \mathbb{E}[\mathcal{S}^{p,\alpha^p}_T]
			= \mathbb{E}\left[\mathcal{W}^{p,\alpha^p}_T+ \int_0^T f_p(s,q_s^{\alpha^p}, \mathbb{E}[q_s^{\alpha^p}], \alpha^p_s,\mathbb{E}[\alpha^p_s];c^{\beta^c})ds\right]\\&
			\overset{\text{i)}}{=} \mathbb{E}\left[g_p(q_T^{\alpha^p}, \mathbb{E}[q_T^{\alpha^p}];c^{\beta^c})+ \int_0^T f_p(s,q_s^{\alpha^p}, \mathbb{E}[q_s^{\alpha^p}], \alpha^p_s,\mathbb{E}[\alpha^p_s];c^{\beta^c})ds\right]
			=\widetilde J_p^\lambda(\alpha^p; c^{\beta^c})\\&
			= \widetilde J_p^\lambda(\alpha^p;\beta^c).
		\end{split}
	\]
	Then, the arbitrariness of $\alpha^p \in \mathcal{A}$ implies that $\mathbb{E}[w_0^p(q_0,\bar q_0)] \geq \sup_{\alpha^p \in \mathcal{A}} \widetilde J_p^\lambda(\alpha^p;c^{\beta^c}) = V_p^\lambda(\beta^c)$.
	
	Performing the same computations with $\alpha^{p,\star}$ instead of $\alpha^p$, by condition iii), we get: $\mathbb{E}[w_0^p(q_0,\bar q_0)] = \widetilde J_p^\lambda(\alpha^{p,\star};\beta^c)$. Then, we have showed that $\alpha^{p,\star}=\mathbf{B}_p(\beta^c)$ is the best response to $\beta^c$. The fact that $\alpha^{c,\star}=\mathbf{B}_c(\beta^p)$ is the best response to $\beta^p$ is proved analogously.
	
	Now, take $\widetilde \alpha^p \in \mathcal{A}$ to be another best response to $\beta^c$. We have
	\[
		\mathbb{E}[\mathcal{S}^{P,\widetilde \alpha^p}_0]
		=\mathbb{E}[w_0^p(q_0,\bar q_0)]
		= V_p^\lambda(\beta^c)
		= \widetilde J_p^\lambda (\widetilde \alpha^p, \beta^c)
		= \mathbb{E}[\mathcal{S}^{P,\widetilde \alpha^p}_T].
	\]
	Then, we conclude that the map $[0,T] \ni t \mapsto \mathbb{E}[\mathcal{S}^{P,\widetilde \alpha^p}_t]$ is constant, since it is nonincreasing and it takes the same value at its extremal points. This reasoning, with a few modifications, can be replicated for $\widetilde \alpha^c$, hence concluding the proof.

\section{Computations of the best response maps}\label{Comp_for_game}

As we have done in Section \ref{Sub-step_1.2} (Sub-step 1.2), we  develop here only the computations for the best response of the producer. The best response of the consumer is obtained following very similar computations.
In this section we show that, setting $w_t^p(q,\bar q)=K_p(t)(q-\bar q)^2+\Lambda_p(t)\bar q^2+2Y^p_t q +R_p(t)$, with $(K_p,\Lambda_p,Y^p,R_p) \in L^{\infty}([0,T], \mathbb R_-)^2\times S^2_{\mathbb F}(\Omega \times [0,T], \mathbb R)\times L^{\infty}([0,T], \mathbb R)$, once $\mathcal{S}^{p,\alpha^p}$ is defined as in the Verification Theorem in Theorem  \ref{Ver_thm}, we have
\[
	\begin{split}
		\frac {d}{dt}\mathbb E[\mathcal{S}^{p,\alpha^p}_t]&
		=\mathbb E \Big[  
		(  {K}_p'(t) + Q_p)(q_t- \mathbb E[q_t])^2 
		+(  \Lambda_p'(t) +Q_p+ \widetilde Q_p)			\mathbb E[q_t]^2
		+2(  {Y^p_t}' + M^{p,c}_t)q_t\\&
		\qquad + {R}_p'(t) + T^{p,c}_t+ \chi^p_t(\alpha^p_t)
		\Big],
	\end{split}
\]
where, for all $t \in [0,T]$, we have set
\begin{equation}
	\left \{ \begin{array}{l}
		\chi^p_t(\alpha^p(t)):= (\alpha^p_t)^\top S_p(t)\alpha^p_t\\
		\qquad +2[U_p(t)(q_t-\mathbb E[q_t])+V_p(t)q_t+\xi^p_t+\bar{\xi}^p_t+O_p(t)]^\top	\alpha^p_t \\
		S_p(t):=N_p+e_2K_p(t)e_2^\top\\
		U_p(t):=K_p(t)e_1\\
		V_p(t):=\Lambda_p(t)e_1\\
		O_p(t):=H_p+ e_1\mathbb{E}[Y^p_t]+e_2\mathbb{E}[Z^{p,W}_t]\\ 
		\xi^p_t:= H_p+ e_1 Y^p_t+e_2 Z^{p,W}_t\\
		\bar{\xi}^p_t:= H_p+ e_1\mathbb{E}[Y^p_t]+e_2\mathbb{E}[Z^{p,W}_t].
	\end{array}\right.
\end{equation}
First of all, we notice 
\begin{equation}
	\frac{d \mathbb{E}[\mathcal{S}^{p,\alpha^p}_t]}{dt}
	= \mathbb{E}\left[ \frac{d}{dt}\mathbb{E}[w_t^p(q_t^{\alpha^p},\mathbb{E}[q_t^{\alpha^p}])]  + f_P(t,q_t^{\alpha^p}, \mathbb{E}[q_t^{\alpha^p}], \alpha_t^p, \mathbb{E}[\alpha_t^p]; c^{\beta^c})\right].
\end{equation}
\\
The dynamics of the state variable controlled by the producer is rewritten as
\begin{align}
	&d\bar q_t^{\alpha^p}= e_1^\top \bar \alpha_t^p dt, \\&
	 d(q_t^{\alpha^p}-\bar q_t^{\alpha^p})= e_1^\top(\alpha_t^p-\bar \alpha_t^p) dt + e_2^\top \alpha_t^p dW_t,
\end{align}
From now on, we write $q_t$ for $q_t^{\alpha^p}$ to simplify the notation. Applying It\^o's formula to $w_t^p(q_t,\mathbb{E}[q_t])$, we get 
\begin{equation*}
	\begin{split}
		dw_t^p(q_t,\mathbb{E}[q_t])&
		= {K}_p'(t)(q_t-\bar{q}_t)^2 dt + K_p(t)[2(q_t-\bar{q}_t)d(q_t-\bar q_t)+ (e_2^\top \alpha_t^p)^2dt] + {\Lambda}_p'(t)(\bar{q}_t)^2 dt\\&
		\quad +2 \Lambda_p(t)\bar{q}_td\bar q_t + 2 q_t dY^p_t +2 Y^p_t dq_t + Z^{p,W}_t e_2^\top \alpha^p_t dt+  {R}_p'(t)dt\\&
		= {K}_p'(t)(q_t-\bar{q}_t)^2 dt + K_p(t)\{2(q_t-\bar{q}_t)[e_1^\top(\alpha^p_t-\bar \alpha^p_t) dt + e_2^\top \alpha_t^p dW_t]+ (e_2^\top \alpha^p_t)^2dt\} \\&
		\quad + {\Lambda}_p'(t)(\bar{q}_t)^2 dt +2 \Lambda_p(t)\bar{q}_t e_1^\top \bar{\alpha}^p_tdt + 2 q_t ( {Y^p_t}' dt + Z^{p,W}_t dW_t +Z_p^B dB_t)\\&
		\quad +2 Y^p_t(e_1^\top \alpha_t^p dt + e_2^\top \alpha_t^p dW_t) + Z^{p,W}_t e_2^\top \alpha^p_t dt + {R}_p'(t)dt\\&
		= [{K}_p'(t)(q_t-\bar{q}_t)^2 + 2K_p(t)(q_t-\bar{q}_t)e_1^\top(\alpha_t^p-\bar{\alpha}_t^p)+K_p(t)(e_2^\top\alpha_t^p)^2+{\Lambda}_p'(t)(\bar{q}_t)^2\\&
		\quad  +2 \Lambda_p(t)\bar{q}_t e_1^\top \bar{\alpha}_t^p + 2 {Y^p_t}' q_t + 2 Y^p_te_1^\top \alpha_t^p +2 Z^{p,W}_t e_2^\top \alpha^p_t + {R}_p'(t) ] dt \\&
		\quad + 2 [ K_p(t)(q_t-\bar{q}_t) e_2^\top \alpha_t^p + Z^{p,W}_t + Y^p_t e_2^\top \alpha_t^p  ] dW_t + 2Z^{p,B}_t  dB_t
	\end{split}
\]
Then, taking its expected value, we obtain
\begin{equation}\label{eq_for_w_game}
	\begin{split}
		\frac{d}{dt}\mathbb{E}[w_t^p(q_t,\bar{q}_t)]&
		=\frac{\mathbb{E}[dw_t^p(q_t,\mathbb{E}[q_t])]}{dt}
		= \mathbb{E}\Big[ {K}_p'(t)(q_t-\bar{q}_t)^2  + 2K_p(t)(q_t-\bar{q}_t)e_1^\top(\alpha_t^p-\bar{\alpha}_t^p)\\&
		\quad +K_p(t)(e_2^\top\alpha_t^p)^2+{\Lambda}_p'(t)(\bar{q}_t)^2  +2 \Lambda_p(t)\bar{q}_t e_1^\top \bar{\alpha}_t^p + 2 {Y^p_t}'q_t + 2 Y^p_te_1^\top \alpha_t^p\\&
		\quad + {R}_p'(t) + 2Z^{p,W}_t e_2^\top \alpha^p_t\Big] \\&
		= \mathbb{E}\Big[{K}_p'(t)(q_t-\bar{q}_t)^2  + 2K_p(t)(q_t-\bar{q}_t)e_1^\top\alpha_t^p +K_p(t)(e_2^\top\alpha_t^p)^2+{\Lambda}_p'(t)(\bar{q}_t)^2\\&
		\quad  +2 \Lambda_p(t)\bar{q}_t e_1^\top {\alpha}_t^p + 2 {Y^p_t}'q_t + 2 Y^p_te_1^\top \alpha_t^p + {R}_p'(t) + 2Z^{p,W}_t e_2^\top \alpha^p_t\Big]\\&
		= \mathbb{E}\Big[{K}_p'(t)(q_t-\bar{q}_t)^2  +{\Lambda}_p'(t)(\bar{q}_t)^2 + 2 {Y^p_t}'q_t  + {R}_p'(t) +K_p(t)(e_2^\top\alpha_t^p)^2 \\&
		\quad + \left[2 (K_p(t)(q_t-\bar{q}_t) + \Lambda_p(t)\bar{q}_t + Y^p_t)e_1 + 2Z^{p,W}_t e_2 \right]^\top \alpha^p_t\Big],
	\end{split}
\end{equation}
where we have used the following semplifications:  $\mathbb{E}[2\Lambda_p(t)e_1^\top \bar \alpha^p_t\bar{q}_t]=\mathbb{E}[2\Lambda_p(t)\bar{q}_t e_1^\top \alpha^p_t]$ and $\mathbb{E}[2K_p(t)e_1^\top \bar \alpha^p_t(q_t - \bar{q}_t)]=2K_p(t)e_1^\top \bar \alpha^p_t\mathbb{E}[q_t-\bar{q}_t]=0$. Moreover, since 
\begin{equation}\label{eq_for_f_game}
	\mathbb{E}[f_p(t,q_t, \bar q_t, \alpha_t^p, \bar \alpha_t^p; c^{\beta^c})]
	= \mathbb{E}[ Q_p (q_t-\bar{q}_t)^2 + (Q_p+\widetilde Q_p)\bar{q}_t^2 +2 M^{p,c}_t q_t+(\alpha^p_t)^\top N_p \alpha^p_t +2H_p^\top \alpha^p_t +T^{p,c}_t],
\end{equation}
by adding up \eqref{eq_for_w_game} and \eqref{eq_for_f_game}, we get
\begin{equation}
	\begin{split}
		\frac{d \mathbb{E}[\mathcal{S}^{p,\alpha^p}_t]}{dt}&
		= \mathbb{E}\left[ \frac{d}{dt}\mathbb{E}[w_t^p(q_t, \bar q_t)]  + f_p(q_t,\bar q_t, \alpha_t^p, \bar \alpha_t^p; c^{\beta^c})\right]\\&
		= \mathbb{E}\Big[ ({K}_p'(t)+  Q_p)(q_t-\bar{q}_t)^2  +({\Lambda}_p'(t) + Q_p+\widetilde Q_p)(\bar{q}_t)^2 + 2 ({Y^p_t}' + M^{p}(c)_t )q_t \\&
		\quad  + {R}_p'(t)  +T^{p}(c)_t  + \chi^p_t(\alpha^p_t)\Big],
	\end{split}
\end{equation}
where we have set
\begin{equation*}
	\begin{split}
		\chi^p_t(\alpha^p_t)
		:&= K_p(t)(e_2^\top \alpha_t^p)^2 + \left\{2 [K_p(t)(q_t-\bar{q}_t) + \Lambda_p(t)\bar{q}_t + Y^p_t] e_1 + 2Z^{p,W}_t e_2 \right\}^\top \alpha^p_t \\&
		\quad  +(\alpha^p_t)^\top N_p \alpha^p_t +2H_p^\top \alpha^p_t\\&
		= \left\{2 [K_p(t)(q_t-\bar{q}_t) + \Lambda_p(t)\bar{q}_t + Y^p_t]\ell_p + 2Z^{p,W}_t e_2 +2H_p\right\}^\top \alpha^p_t\\&
		\quad +(\alpha^p_t)^\top (N_p + e_2 K_p(t)e_2^\top) \alpha^p_t\\&
		= 2[U_p(t)(q_t-\mathbb E[q_t])+V_p(t)q_t+\xi^p_t+\bar{\xi}^p_t+O_p(t)]^\top	\alpha^p_t\\&
		\quad + (\alpha^p_t)^\top S_p(t)\alpha^p_t,
	\end{split}
\end{equation*}
with
\begin{equation*}
	\left \{ \begin{array}{l}
		S_p(t):=N_p+e_2K_p(t)e_2^\top\\
		U_p(t):=K_p(t)e_1\\
		V_p(t):=\Lambda_p(t)e_1\\
		O_p(t):=H_p+ e_1\mathbb{E}[Y^p_t]+e_2\mathbb{E}[Z^{p,W}_t]\\ 
		\xi^p_t:= H_p+ e_1 Y^p_t+e_2 Z^{p,W}_t\\
		\bar{\xi}^p_t:= H_p+ e_1\mathbb{E}[Y^p_t]+e_2\mathbb{E}[Z^{p,W}_t].
	\end{array}\right.
\end{equation*}

\section{Computations of the equilibrium payoffs}\label{app:valfunc}

In this section we perform some computations to get a more explicit formula for the objective functionals at the equilibrium in Theorem \ref{thm_Nash}.
In particular, we find explicit expressions for $R_p(0)$ and $R_c(0)$.
In all the following computations we are using the optimal strategies but we are suppressing the stars in the notation for the sake readability (e.g. we write $u_t$ instead of $u^* _t$ and so on). For the same reason we are suppressing the dependency on time when clear  from the context.

\begin{proposition}
It holds that
\begin{align*}
R^{(\lambda)}_p(0)=& \int_0^T \left[ \frac{2}{k_p} \E [(Y^p_u)^2]   -\eta_p \lambda^2 \gamma^2\rho_c^2  \V [c_u] 
 +\frac{2\big(     \pi_{11}(u) z_u  + \frac{\ell_p \sigma_p}{2}\big )^2}{\ell_p-2K_p(u)}  \right] du - \lambda \gamma \rho_c  \bar{c}_T, \\
R^{(\lambda)}_c(0)= &   \int_0^T \left[   \frac{2}{k_c}  \E [(Y^c_u)^2]  
-\eta_c \lambda^2 \rho_p^2 \V [q_u]
 + \frac{2\big( \pi_{22}(u) y_u + \frac{\ell_c \sigma_c}{2}\big)^2}{\ell_c-2K_c(u)}   \right ]du - \lambda \rho_p \bar{q}_T,
\end{align*}	
where 
\begin{align*}
	&d\bar c_t = \frac{2}{k_c}  \Big[ \big( \Lambda_c +\widehat \pi_{22} \big) \bar c_t  + \widehat \pi_{21}   \bar q_t +  h_2   \Big], \\
	&d\bar q_t = \frac{2}{k_p} \Big[\widehat \pi_{12}  \bar c_t + \big( \Lambda_p +\widehat \pi_{11} \big) \bar q_t  +   h_1 \Big], \\
	&d\E[c^2_t]  = \frac{4}{k_c} \Big[ \left( K_c  + \pi_{22}  \right) \big( \E[c^2_t] - \bar c_t^2 \big)  +   \pi_{21}  \big( \E[ c_t q_t] - \bar c_t \bar q_t \big)  +  \big[ \Lambda_c  + \widehat \pi_{22} \big] \bar c^2_t 
 +  \widehat \pi_{21} \bar c_t \bar q_t +  h_2 \bar c_t \Big] dt+   y_t^2   dt, \\
	&d\E[q^2_t]  = \frac{4}{k_p} \Big[\left( K_p  + \pi_{11}  \right) \big( \E[q^2_t] - \bar q_t^2 \big) + \pi_{12}  \big( \E[ c_t q_t] - \bar c_t \bar q_t \big)
+ \big( \Lambda_p  + \widehat \pi_{11}  \big) \bar q^2_t + \widehat \pi_{12}  \bar c_t \bar q_t + h_1 \bar q_t \Big] dt+   z_t^2   dt, \\
	& d\E[ c_t q_t] =   \frac{2}{k_p} \Big[ \big( K_p  + \pi_{11}  \big) \big( \E[ c_t q_t] - \bar c_t \bar q_t \big) 
+ \pi_{12} \big( \E[c^2_t] - \bar c^2_t \big) + \big( \Lambda_p + \widehat \pi_{11}   \big) \bar c_t \bar q_t  + \widehat  \pi_{12}   \bar c^2_t + h_1  \bar c_t \Big]dt\\
& \hspace{1.7cm} +  \frac{2}{k_c} \Big[ \big(  K_c + \pi_{22} \big) \big( \E[ c_t q_t] - \bar c_t \bar q_t \big) 
+ \pi_{21} \big( \E[q^2_t] - \bar q^2_t \big) + \big( \Lambda_c + \widehat \pi_{22} \big) \bar c_t \bar q_t  + \widehat  \pi_{21}  \bar q^2_t + h_2  \bar q_t \Big]dt,
\end{align*}
\begin{align*}
	&  d\E[(Y^p_t)^2] = - 2 \Big\{   
\frac12 s_0 \bar Y^p_t + \frac12 \gamma \rho_c \E[Y^p_t c_t]
+ \rho_p \gamma \rho_c \eta_p \lambda^2 \big( \E[Y^p_t c_t ] - \bar Y^p_t \bar c_t \big)
+ \frac{2}{k_p} \left[ K_p  \E\big[ (Y^p_t)^2 - (\bar Y^p_t)^2 \big] + \Lambda_p (\bar Y^p_t)^2 \right]
\Big\} dt \\
& \hspace{2.06cm} +  \big(\pi_{11}^2 z^2_t + \pi_{12}^2 y^2_t \big) dt,\\
	&  d\E[(Y^c_t)^2] = - 2 \Big\{ \frac{p_0+p_1 s_0 -\gamma(s_0+\delta)}{2} \bar Y^c_t + \frac{\rho_p(\gamma-p_1)}{2}\E[Y^c_t q_t]
+ \rho_p \gamma \rho_c \eta_c \lambda^2 \big( \E[Y^c_t q_t ] - \bar Y^c_t \bar q_t \big)\\
& \hspace{2.05cm} + \frac{2}{k_c} \left[  K_c  \E\big[(Y^c_t)^2 - (\bar Y^c_t)^2 \big] + \Lambda_c (\bar Y^c_t)^2 \right] \Big\} dt +  \big(\pi_{21}^2 z^2_t + \pi_{22}^2 y^2_t \big) dt,
\end{align*}
with
\begin{align*}
	& \bar Y^p_t = \widehat \pi_{11} \bar q_t + \widehat \pi_{12} \bar c_t +  h_1,\quad
	 \bar Y^c_t = \widehat \pi_{21} \bar q_t + \widehat \pi_{22} \bar c_t +  h_2,\\
	& \E[Y^p_t c_t ] = \pi_{11} \big( \E[c_t q_t] - \bar c_t \bar q_t \big) + \pi_{12} \big( \E[ c^2_t] - \bar c^2_t \big) + \widehat \pi_{11} \bar c_t \bar q_t + \widehat \pi_{12} \bar c_t^2 + h_1 \bar c_t, \\
	& \E[Y^c_t q_t ] = \pi_{21} \big( \E[q_t^2] - \bar q_t^2 \big) + \pi_{22} \big( \E[ q_t c_t] - \bar c_t \bar{q}_t \big) + \widehat \pi_{21} \bar q_t^2 + \widehat \pi_{22} \bar c_t \bar q_t + h_2 \bar q_t,
\end{align*}
and terminal conditions
\begin{align*}
	(Y^p_T)^2 = \frac{\lambda^2 \rho_p^2 }{4},
	\quad 
	(Y^c_T)^2 = \frac{\lambda^2 \gamma^2\rho_c^2 }{4}.
\end{align*}
\end{proposition}

\begin{proof}

\noindent For the terms  $\bar q_t = \mathbb{E}[q_t]$ and $\bar c_t =$ $\E[c_T]$ we have
\begin{align*}
&d\bar q_t = \bar u_t dt =  \frac{2}{k_p} \Big[ \big( \Lambda_p +\widehat \pi_{11} \big) \bar q_t  + \widehat \pi_{12}  \bar c_t +  h_1 \Big] dt, \quad
d\bar c_t = \bar v_t dt =  \frac{2}{k_c}  \Big[ \big( \Lambda_c +\widehat \pi_{22} \big) \bar c_t  + \widehat \pi_{21}   \bar q_t +  h_2   \Big] dt, 
\end{align*}
so we have a 2-dimensional ODE giving $\bar c_T$ and $\bar q_T$. 

For the terms   $\V[q_t]$ and $\V[c_t]$, we have
\begin{align*}
&\V[c_t] = \E[c^2_t] - \bar c_t^2, \quad
	d\E[c^2_t] 
	= \big( 2 \E[c_t v_t] + (y_t)^2 \big) dt,\\
&  \V[q_t] = \E[q^2_t] - \bar q_t^2, \quad
 d\E[q^2_t] 
 = \big( 2 \E[q_t u_t] + (z_t)^2 \big) dt,
\end{align*}
because $z_t$ and $y_t$ are deterministic. Further, 
\begin{align*}
&\E[c_t v_t] = \frac{2}{k_c} \Big[\left( K_c  + \pi_{22}  \right) \big( \E[c^2_t] - \bar c_t^2 \big) + \pi_{21}  \big( \E[ c_t q_t] - \bar c_t \bar q_t \big)
+ \big( \Lambda_c  + \widehat \pi_{22}  \big) \bar c^2_t + \widehat \pi_{21}  \bar c_t \bar q_t + h_2 \bar c_t \Big], \\
&d\E[c^2_t]  = \frac{4}{k_c} \Big[ \left( K_c  + \pi_{22}  \right) \big( \E[c^2_t] - \bar c_t^2 \big)  +   \pi_{21}  \big( \E[ c_t q_t] - \bar c_t \bar q_t \big)  +  \big[ \Lambda_c  + \widehat \pi_{22} \big] \bar c^2_t 
 +  \widehat \pi_{21} \bar c_t \bar q_t +  h_2 \bar c_t \Big] dt+   y_t ^2   dt, \\
&\E[q_t u_t] = \frac{2}{k_p} \Big[\left( K_p  + \pi_{11}  \right) \big( \E[q^2_t] - \bar q_t^2 \big) + \pi_{12}  \big( \E[ c_t q_t] - \bar c_t \bar q_t \big)
+ \big( \Lambda_p  + \widehat \pi_{11}  \big) \bar q^2_t + \widehat \pi_{12}  \bar c_t \bar q_t + h_1 \bar q_t \Big], \\
&d\E[q^2_t]  = \frac{4}{k_p} \Big[\left( K_p  + \pi_{11}  \right) \big( \E[q^2_t] - \bar q_t^2 \big) + \pi_{12}  \big( \E[ c_t q_t] - \bar c_t \bar q_t \big)
+ \big( \Lambda_p  + \widehat \pi_{11}  \big) \bar q^2_t + \widehat \pi_{12}  \bar c_t \bar q_t + h_1 \bar q_t \Big] dt+   z_t^2   dt, 
\end{align*}
and we have for $\E[c_t q_t]$, that $d\E[ c_t q_t] = \E \big( c_t u_t   + q_t v_t \big) dt$, so that
\begin{align*}
 d\E[ c_t q_t] =  & \frac{2}{k_p} \Big[ \big( K_p  + \pi_{11}  \big) \big( \E[ c_t q_t] - \bar c_t \bar q_t \big) 
+ \pi_{12} \big( \E[c^2_t] - \bar c^2_t \big) + \big( \Lambda_p + \widehat \pi_{11}   \big) \bar c_t \bar q_t  + \widehat  \pi_{12}   \bar c^2_t + h_1  \bar c_t \Big] \\
+ &  \frac{2}{k_c} \Big[ \big(  K_c + \pi_{22} \big) \big( \E[ c_t q_t] - \bar c_t \bar q_t \big) 
+ \pi_{21} \big( \E[q^2_t] - \bar q^2_t \big) + \big( \Lambda_c + \widehat \pi_{22} \big) \bar c_t \bar q_t  + \widehat  \pi_{21}  \bar q^2_t + h_2  \bar q_t \Big]. 
\end{align*}
For the term $\V[Y^p_t]$, we have $\V[Y^p_t] + \E[Y^p_t]^2 = \E[(Y^p_t)^2] $, where
\begin{align*}
&  Y^p_t = \pi_{11} (q_t - \bar q_t) + \pi_{12} (c_t - \bar c_t) + \widehat \pi_{11} \bar q_t + \widehat \pi_{12} \bar c_t + h_1(t), \quad \bar Y^p_t = \widehat \pi_{11} \bar q_t + \widehat \pi_{12} \bar c_t +  h_1, \\
&  \E[Y^p_t c_t ] = \pi_{11} \big( \E[c_t q_t] - \bar c_t \bar q_t \big) + \pi_{12} \big( \E[ c^2_t] - \bar c^2_t \big) + \widehat \pi_{11} \bar c_t \bar q_t + \widehat \pi_{12} \bar c_t^2 + h_1 \bar c_t, \\
&  d\E[(Y^p_t)^2] = - 2 \Big\{   
\frac12 s_0 \bar Y^p_t + \frac12 \gamma \rho_c \E[Y^p_t c_t]
+ \rho_p \gamma \rho_c \eta_p \lambda^2 \big( \E[Y^p_t c_t ] - \bar Y^p_t \bar c_t \big)
+ \frac{2}{k_p} \left[ K_p  \E\big[ (Y^p_t)^2 - (\bar Y^p_t)^2 \big] + \Lambda_p (\bar Y^p_t)^2 \right]
\Big\} dt \\
& \hspace{2.1cm}+  \big(\pi_{11}^2 z^2_t + \pi_{12}^2 y^2_t \big) dt,\quad  (Y^p_T)^2 = \frac14 \lambda^2 \rho_p^2, 
\end{align*}
where we have exploited the representation of $Y^p$ in Equation  \eqref{Sys_p}. Analogously, exploiting the representation of $Y^c$ in Equation \eqref{Sys_c}, we get
\begin{align*}
&  Y^c_t = \pi_{21} (q_t - \bar q_t) + \pi_{22} (c_t - \bar c_t) + \widehat \pi_{21} \bar q_t + \widehat \pi_{22} \bar c_t + h_2(t), \quad \bar Y^c_t = \widehat \pi_{21} \bar q_t + \widehat \pi_{22} \bar c_t +  h_2, \\
& + \frac{2}{k_c} \left[  K_c  \E\big[(Y^c_t)^2 - (\bar Y^c_t)^2 \big] + \Lambda_c (\bar Y^c_t)^2 \right] \Big\} dt, \\
&  \E[Y^c_t q_t ] = \pi_{21} \big( \E[q_t^2] - \bar q_t^2 \big) + \pi_{22} \big( \E[ q_t c_t] - \bar c_t \bar{q}_t \big) + \widehat \pi_{21} \bar q_t^2 + \widehat \pi_{22} \bar c_t \bar q_t + h_2 \bar q_t, \\
&  d\E[(Y^c_t)^2] = - 2 \Big\{ \frac{p_0+p_1 s_0 -\gamma(s_0+\delta)}{2} \bar Y^c_t + \frac{\rho_p(\gamma-p_1)}{2}\E[Y^c_t q_t]
+ \rho_p \gamma \rho_c \eta_c \lambda^2 \big( \E[Y^c_t q_t ] - \bar Y^c_t \bar q_t \big)\\
& \hspace{2.1cm} + \frac{2}{k_c} \left[  K_c  \E\big[(Y^c_t)^2 - (\bar Y^c_t)^2 \big] + \Lambda_c (\bar Y^c_t)^2 \right] \Big\} dt +  \big(\pi_{21}^2 z^2_t + \pi_{22}^2 y^2_t \big) dt,\quad  (Y^p_T)^2 = \frac14 \lambda^2 \gamma^2\rho_c^2, 
\end{align*}
 Summing up, we have obtained a backward ODE for $\E[(Y^p_t)^2]$ and $\E[(Y^c_t)^2]$. Finally, we have  
\begin{align*}
& Z^{p,W}_t = \pi_{11} z_t, \quad \V[Z^{p,W}_t] = 0, \quad \E[Z^{p,W}_t] =  \pi_{11} z_t
\end{align*}
and
\begin{align*}
& Z^{c,B}_t = \pi_{22} y_t, \quad \V[Z^{c,B}_t] = 0, \quad \E[Z^{c,B}_t] =  \pi_{22} y_t.
\end{align*}
Recalling that
\begin{align*}
R_p^{(\lambda)}(t)=& - \lambda \gamma \rho_c \mathbb E [c_T]+\int_t^T \Bigg[-\eta_p \lambda^2 \gamma^2\rho_c^2 \mathbb V [c_u] +\frac{2}{k_p}\left(\mathbb V [Y^p_u]+\mathbb E [Y^p_u]^2\right)\\&
+\frac{2}{\ell_p-2K_p(u)}\left(\mathbb V[Z^{p,W}_u]+\left(\mathbb E [Z^{p,W}_u]+\frac{\ell_p \sigma_p}{2}\right)^2\right)\Bigg]du,
\end{align*}
and analogously
\begin{align*}
R_c^{(\lambda)}(t)=& - \lambda \rho_p \mathbb E [q_T]+\int_t^T \Bigg[-\eta_c \lambda^2 \rho_p^2 \mathbb V [q_u] +\frac{2}{k_c}\left(\mathbb V [Y^c_u]+\mathbb E [Y^c_u]^2\right)\\&
+\frac{2}{\ell_c-2K_c(u)}\left(\mathbb V[Z^{c,B}_u]+\left(\mathbb E [Z^{c,B}_u]+\frac{\ell_c \sigma_c}{2}\right)^2\right)\Bigg]du,
\end{align*}
the results follow.
\end{proof}


\bibliographystyle{plain}

\end{document}